\documentclass[11pt,leqno]{article}
\normalfont
\renewcommand{\baselinestretch}{1.25}

\usepackage{graphicx}
\usepackage[dvipsnames]{xcolor}

\newcommand{\IntroSection}{1}
\newcommand{\IntroFibFig}{Figure 1.1}
\newcommand{\TheoremList}{Theorem 4.2/Corollary 4.3/Theorem 5.1}
\newcommand{\ArraySection}{2}
\newcommand{\RankProposition}{Proposition 2.1}
\newcommand{\EnumerativeTheorem}{Theorem 2.2}
\newcommand{\TableFig}{Figure 2.1}
\newcommand{\EnumerativeCorollary}{Corollary 2.3}
\newcommand{\SetupSection}{3}
\newcommand{\SetupProposition}{Proposition 3.1}
\newcommand{\SolitaryLatticeProposition}{Proposition 3.2}
\newcommand{\MainResultSection}{4}
\newcommand{\DDTheorem}{Theorem 4.1}
\newcommand{\RibbonFigure}{Figure 4.1}
\newcommand{\GTFig}{Figure 4.2}
\newcommand{\FibGTTheorem}{Theorem 4.2}
\newcommand{\FibGTCorollary}{Corollary 4.3}
\newcommand{\SCDTheorem}{Theorem 4.4}
\newcommand{\RGFSection}{5}
\newcommand{\RGFTheorem}{Theorem 5.1}
\newcommand{\EndSection}{6}
\newcommand{\LatticeFigList}{Figures 1.1 and 4.1}
\newcommand{\FibSolitaryFig}{Figure 6.1}
\newcommand{\FibSpineFig}{Figure 6.2}
\newcommand{\FibSolitaryProposition}{Proposition 6.1}
\newcommand{\FibSolitaryOpen}{Open Problem 6.2}
\newcommand{\FibMotzkinFig}{Figure 6.3}
\newcommand{\MotzkinConjecture}{Conjecture 6.3}
\newcommand{\MotzkinOpen}{Open Problems 6.4}

\usepackage{amsfonts}
\usepackage{mathrsfs}
\usepackage{rotating}
\usepackage{mathdots}

\newfont{\mysmallscbolditalics}{ecoc0500 at 10pt}
\newcommand{\mymathfrak}[1]{\mbox{\mysmallscbolditalics #1}}

\newfont{\myscbolditalics}{ecoc0500 at 11pt}

\newfont{\mybolditalics}{ecbi0500 at 11pt}

\newcommand{\myqx}{\mbox{\mybolditalics x}}
\newcommand{\myqy}{\mbox{\mybolditalics y}}
\newcommand{\myqz}{\mbox{\mybolditalics z}}
\newcommand{\myqX}{\mbox{\mybolditalics X}}
\newcommand{\myqY}{\mbox{\mybolditalics Y}}
\newcommand{\myqh}{\mbox{\mybolditalics h}}
\newcommand{\myqP}{\mbox{\mybolditalics P}}

\newfont{\eulercursive}{eurm10 at 11pt}

\newcommand{\mya}{\mbox{\eulercursive a}}
\newcommand{\mys}{\mbox{\eulercursive s}}
\newcommand{\mym}{\mbox{\eulercursive m}}

\newfont{\smalleulercursive}{eurm10 at 9pt}

\newfont{\smallereulercursive}{eurm10 at 7pt}

\newcommand{\mysmallerm}{\mbox{\smallereulercursive m}}

\newfont{\myslantcyrillic}{wncyi10 at 11pt}

\newcommand{\QED}{\raisebox{0.5mm}{\fbox{\rule{0mm}{1.5mm}\ }}}

\newcounter{myfn}[page]
\renewcommand{\thefootnote}{\fnsymbol{footnote}}

\newcounter{rone}
\newcounter{rtwo}
\newcounter{rthree}
\newcounter{rfour}
\newcounter{rfive}
\newcounter{rsix}
\newcounter{rseven}
\newcommand{\myA}{\mbox{\sffamily A}}

\newcommand{\mytinyA}{\mbox{\tiny \sffamily A}}
\newcommand{\myB}{\mbox{\sffamily B}}

\newcommand{\myC}{\mbox{\sffamily C}}

\newcommand{\mysmallP}{\mbox{\footnotesize \sffamily P}}

\newcommand{\mytinyP}{\mbox{\tiny \sffamily P}}

\newcommand{\mysmallQ}{\mbox{\footnotesize \sffamily Q}}

\newcommand{\mytinyQ}{\mbox{\tiny \sffamily Q}}

\newcommand{\mysmallR}{\mbox{\footnotesize \sffamily R}}

\newcommand{\myvarZ}{\mbox{\scriptsize \sffamily Z}}
\newcommand{\mypart}{\mbox{\small \sffamily Part}}
\newcommand{\mysmallerpart}{\mbox{\tiny \sffamily Part}}

 \newcommand{\relt}{\mathbf{r}}
\newcommand{\selt}{\mathbf{s}} \newcommand{\telt}{\mathbf{t}}
\newcommand{\uelt}{\mathbf{u}} \newcommand{\velt}{\mathbf{v}}

\newcommand{\wt}{\mbox{\sffamily wt}}
\newcommand{\smallwt}{\mbox{\scriptsize \sffamily wt}}

\newcommand{\WGF}{\mbox{\sffamily WGF}}
\newcommand{\RGF}{\mbox{\sffamily RGF}}

\newcommand{\altFibpoly}[1]{{\overline{F}_{#1}}}

\newcommand{\comp}{\mbox{\sffamily comp}}

\newcommand{\myarrow}[1]{\stackrel{#1}{\rightarrow}}

\newcommand{\mylongarrow}[1]{\stackrel{#1}{\longrightarrow}}

\newcommand{\NEEdgeLabelForLatticeI}[1]{
\setlength{\unitlength}{1.5cm}
\begin{picture}(0,0)
\put(-0.25,0){
\begin{picture}(0,0)
\put(0.4,0.4){\footnotesize #1} 
\end{picture}
}
\end{picture}
}

\newcommand{\NWEdgeLabelForLatticeI}[1]{
\setlength{\unitlength}{1.5cm}
\begin{picture}(0,0)
\put(-0.25,0){
\begin{picture}(0,0)
\put(-0.525,0.4){\footnotesize #1} 
\end{picture}
}
\end{picture}
}

\newcommand{\VerticalEdgeLabelForLatticeI}[1]{
\setlength{\unitlength}{1.5cm}
\begin{picture}(0,0)
\put(-0.25,0){
\begin{picture}(0,0)
\put(-0.05,0.4){\footnotesize #1} 
\end{picture}
}
\end{picture}
}

\newcommand{\VertexkTupleFib}[5]{
\setlength{\unitlength}{1.5cm}
\begin{picture}(0,0)
\put(-0.25,0){
\begin{picture}(0,0)
\put(0,0){\circle*{0.1}} 
\put(#4,#5){\scriptsize $(#1,#2,#3)$}
\end{picture}
}
\end{picture}
}

\newcommand{\VertexBlankFib}{
\setlength{\unitlength}{1.5cm}
\begin{picture}(0,0)
\put(-0.25,0){
\begin{picture}(0,0)
\put(0,0){\circle*{0.1}} 
\end{picture}
}
\end{picture}
}

\newcommand{\VertexTableauFibTwo}[6]{
\setlength{\unitlength}{1.5cm}
\begin{picture}(0,0)
\put(-0.25,0){
\begin{picture}(0,0)
\put(0,0){\circle*{0.1}} 
\put(#5,#6){\setlength{\unitlength}{0.25cm}\begin{picture}(0,0)\put(0,0){\line(0,1){1}} \put(1,0){\line(0,1){2}} \put(2,0){\line(0,1){2}} \put(3,1){\line(0,1){1}} \put(0,0){\line(1,0){2}} \put(0,1){\line(1,0){3}} \put(1,2){\line(1,0){2}} \put(0.25,0.25){\tiny #1} \put(1.25,1.25){\tiny #2} \put(1.25,0.25){\tiny #3} \put(2.25,1.25){\tiny #4} \end{picture}}
\end{picture}
}
\end{picture}
}

\begin{document}
\pagenumbering{arabic}
\thispagestyle{empty}%
\vspace*{-0.7in}
\noindent
{\scriptsize To appear in {\em Involve, a Journal of Mathematics}}\hfill {\scriptsize May 15, 2022} 

\begin{center}
{\large \bf Symmetric Fibonaccian distributive lattices and\\ 
representations of the special linear Lie algebras} 

\vspace*{0.05in}
\renewcommand{\thefootnote}{1}
Robert G.\ Donnelly,\footnote{Department of Mathematics and Statistics, Murray State
University, Murray, KY 42071\\ 
\hspace*{0.25in}Email: {\tt rob.donnelly@murraystate.edu}}  
\renewcommand{\thefootnote}{2} 
\hspace*{-0.07in}Molly W.\ Dunkum,\footnote{Department of Mathematics, Western Kentucky University, Bowling Green, KY 42101\\ 
\hspace*{0.25in}Email: {\tt molly.dunkum@wku.edu}} 
\renewcommand{\thefootnote}{3} 
\hspace*{-0.07in}Sasha V.\ Malone,\footnote{Department of Mathematics, Western Kentucky University, Bowling Green, KY 42101\\ 
\hspace*{0.25in}Email: {\tt sverona@sverona.xyz}}
\renewcommand{\thefootnote}{4} 
\hspace*{-0.07in}and Alexandra Nance\footnote{Department of Mathematics and Statistics, Murray State
University, Murray, KY 42071\\ 
\hspace*{0.25in}Email: {\tt alexandrannance@gmail.com}}
\end{center} 

\begin{abstract}
We present a family of rank symmetric diamond-colored distributive lattices that are naturally related to the Fibonacci sequence and certain of its generalizations. 
These lattices re-interpret and unify descriptions of some un- or differently-colored lattices found variously in the literature. 
We demonstrate that our symmetric Fibonaccian lattices naturally realize certain (often reducible) representations of the special linear Lie algebras, with weight basis vectors realized as lattice elements and Lie algebra generators acting along the covering digraph edges of each lattice. 
We present evidence that each such weight basis possesses certain distinctive extremal properties.  
We provide new descriptions of the lattice cardinalities and rank generating functions and offer several conjectures/open problems. 
Throughout, we make connections with integer sequences from the OEIS.

\begin{center}
\renewcommand{\thefootnote}{$\mbox{\ }$} 
{\small \bf Mathematics Subject Classification:}\ {\small 05E15 
(20F55, 17B10)}\\
{\small \bf Keywords:}\ diamond-colored distributive lattice, rank generating function, skew-shaped semistandard tableau, skew Schur function, skew-tabular lattice, special linear Lie algebra representation, weight basis supporting graph / representation diagram\footnote{\ \ }  

\end{center} 
\end{abstract}

\setcounter{footnote}{\value{rone}}
\noindent {\bf \S \IntroSection\ Introduction}. 
Many rich interactions between combinatorics and representation theory have been discovered over the past century or so and continue to be productively explored. 
One simple example is that many prominent integer patterns naturally occur as dimensions of (interesting families of) Lie algebra representations.  
Notably, the $n^{\mbox{\tiny th}}$ power of two (OEIS sequence A000079 \cite{OEIS}\footnote{For the remainder of the paper, `\cite{OEIS}' is the reference for all OEIS-identified sequences.}) is the dimension of the minuscule representation of the type $\myB_{n}$ odd orthogonal Lie algebra $\mathfrak{so}(2n+1,\mathbb{C})$; the $(n+1)^{\mbox{\tiny st}}$ Catalan number $\frac{1}{n+2}{2(n+1) \choose n+1}$ (OEIS-A000108) is the dimension of the $n^{\mbox{\tiny th}}$ (non-minuscule) fundamental representation of the type $\myC_{n}$ symplectic Lie algebra $\mathfrak{sp}(2n,\mathbb{C})$; and the binomial coefficient ${n \choose k}$ (OEIS-A007318) occurs as the dimension of the $k^{\mbox{\tiny th}}$ fundamental representation of the type $\myA_{n-1}$ special linear Lie algebra $\mathfrak{sl}(n,\mathbb{C})$. 
However, amongst such occurrences of representation dimensions naturally enumerated by interesting integer patterns, the Fibonacci numbers (OEIS-A000045\footnote{Here, we regard the Fibonacci sequence to be $f_{-1} = 0, f_{0}=1$, $f_{1}=1$, $f_{2}=2$, $f_{3}=3$, $f_{4}=5$, etc.}) are noticeably absent. 
A principal aim of this paper is to (partially) fill this gap.

\newpage
To do so, we produce some finite diamond-colored distributive lattices\footnote{We follow the conventions of \cite{DonDiamond} for order-theoretic concepts relating to diamond-colored distributive lattices.} that re-interpret certain finite distributive lattices appearing in many different contexts in the literature dating back to the 1970's.  
We define these lattices as follows. 
Fix positive integers $n$ and $k$, with $n \geq 2$. Let 
\[L^{\mbox{\tiny Fib}}(n,k) := \left\{T=(T_{1},\ldots,T_{k})\, \rule[-5mm]{0.2mm}{12mm}\, \begin{array}{c}T_{j} \in \{\mbox{\small $(j-1)n+1,(j-1)n+2,\ldots,jn$}\}\ \mbox{for all}\ j \in \{1,\ldots,k\}\\ \mbox{and}\ T_{j+1} \ne T_{j}+1\ \mbox{for all}\ j \in \{1,\ldots,k-1\}\end{array}\right\},\]
a set of positive integer $k$-tuples satisfying certain conditions. 
Order these $k$-tuples by {\em reverse component-wise comparison} such that $S \leq T$ in $L^{\mbox{\tiny Fib}}(n,k)$ for $S = (S_{1},\ldots,S_{k})$ and $T=(T_{1},\ldots,T_{k})$ if and only if $S_{j} \geq T_{j}$ for any $j \in \{1,\ldots,k\}$. 
Observe that $T$ covers $S$ in the resulting covering digraph (aka Hasse diagram) if and only if there exists some $l \in \{1,\ldots,k\}$ such that $S_{l} = T_{l}+1$ while $S_{j} = T_{j}$ for all $j \ne l$. 
In this case, we view the covering relation as a directed edge $S \myarrow{i} T$ `colored' by $i$, where $i = T_{l}\ \mbox{\small mod}\, n$ if $l$ is odd and $i = n - (T_{l}\ \mbox{\small mod}\, n)$ if $l$ is even. 
It is easy to see that $L^{\mbox{\tiny Fib}}(n,k)$ is a distributive lattice such that on any `diamond' of edges \parbox{1.4cm}{\begin{center}
\setlength{\unitlength}{0.2cm}
\begin{picture}(5,3)
\put(2.5,0){\circle*{0.5}} \put(0.5,2){\circle*{0.5}}
\put(2.5,4){\circle*{0.5}} \put(4.5,2){\circle*{0.5}}
\put(0.5,2){\line(1,1){2}} \put(2.5,0){\line(-1,1){2}}
\put(4.5,2){\line(-1,1){2}} \put(2.5,0){\line(1,1){2}}
\put(1.25,0.55){\em \small k} \put(3.2,0.7){\em \small l}
\put(1.2,2.7){\em \small i} \put(3.25,2.55){\em \small j}
\put(3,-0.75){\footnotesize $\relt$} \put(5.25,1.75){\footnotesize $\telt$}
\put(3,4){\footnotesize $\uelt$} \put(-1,1.75){\footnotesize $\selt$}
\end{picture} \end{center}} in the covering digraph, necessarily $k=j$ and $l=i$.
That is, $L^{\mbox{\tiny Fib}}(n,k)$ is a diamond-colored distributive lattice.  
Also, $L^{\mbox{\tiny Fib}}(n,k)$ is self-dual (and therefore rank symmetric) when we disregard edge colors. 
See \IntroFibFig\ for an example. 

We call $L^{\mbox{\tiny Fib}}(n,k)$ a {\em symmetric Fibonaccian (distributive) lattice} (or {\em SFL} for short) and refer to its elements as $(n,k)$-{\em Fibonaccian strings} or simply {\em Fibonaccian strings} when the context is clear. 
Throughout, we regard $L^{\mbox{\tiny Fib}}(n,0)$ to be a one-element poset and $L^{\mbox{\tiny Fib}}(n,-1)$ to be empty. 
When $n=2$ and $k \geq 0$, notice that each $L^{\mbox{\tiny Fib}}(2,k)$ is a chain with $k+1$ elements. 
When $n=3$, we will show that the sequence of lattice sizes, with starting index $k=0$, is $1, 3, 8, 21, 55, 144, \ldots$, coinciding with the Fibonacci subsequence $\{f_{2m+1}\}_{m \geq 0}$ (cf.\ OEIS-A001906).  
We call numbers in this latter subsequence the {\em symmetric Fibonacci numbers}, as they naturally enumerate the sizes of our SFLs when $n=3$. 
For each $n > 3$, we regard the sequence $\left|\rule[-1.5mm]{0mm}{4.75mm}L^{\mbox{\tiny Fib}}(n,1)\right|, \left|\rule[-1.5mm]{0mm}{4.75mm}L^{\mbox{\tiny Fib}}(n,2)\right|, \left|\rule[-1.5mm]{0mm}{4.75mm}L^{\mbox{\tiny Fib}}(n,3)\right|, \ldots$ to be a generalization of the sequence of symmetric Fibonacci numbers. 

\begin{figure}[t]
\begin{center}
{\small {\bf \IntroFibFig:}  The covering digraph (aka Hasse diagram) for the symmetric Fibonaccian lattice $L^{\mbox{\tiny Fib}}(3,3)$.}

\setlength{\unitlength}{1.5cm}
\begin{picture}(4,6.5)
\put(0,0){\qbezier(2,6)(3.25,5)(4.5,4)}
\put(0,0){\qbezier(2,6)(1.5,5.5)(1,5)}
\put(0,0){\qbezier(2,6)(2,5.5)(2,5)}
\put(0,0){\qbezier(3.25,5)(2.75,4.5)(2.25,4)}
\put(0,0){\qbezier(3.25,5)(3.25,4)(3.25,3)}
\put(0,0){\qbezier(4.5,4)(4,3.5)(3.5,3)}
\put(0,0){\qbezier(4.5,4)(4.5,3)(4.5,2)}
\put(0,0){\qbezier(1,5)(2.25,4)(3.5,3)}
\put(0,0){\qbezier(1,5)(1,4.5)(1,4)}
\put(0,0){\qbezier(2,5)(3.25,4)(4.5,3)}
\put(0,0){\qbezier(2,5)(1,4)(0,3)}
\put(0,0){\qbezier(1,4)(2.25,3)(3.5,2)}
\put(0,0){\qbezier(2.25,4)(2.25,3)(2.25,2)}
\put(0,0){\qbezier(3.5,3)(3.5,2)(3.5,1)}
\put(0,0){\qbezier(3.25,4)(2.25,3)(1.25,2)}
\put(0,0){\qbezier(4.5,3)(3.5,2)(2.5,1)}
\put(0,0){\qbezier(0,3)(1.25,2)(2.5,1)}
\put(0,0){\qbezier(1.25,2)(1.25,1.5)(1.25,1)}
\put(0,0){\qbezier(2.5,1)(2.5,0.5)(2.5,0)}
\put(0,0){\qbezier(3.25,3)(2.25,2)(1.25,1)}
\put(0,0){\qbezier(3.25,3)(3.875,2.5)(4.5,2)}
\put(0,0){\qbezier(4.5,2)(3.5,1)(2.5,0)}
\put(0,0){\qbezier(2.25,2)(2.875,1.5)(3.5,1)}
\put(0,0){\qbezier(1.25,1)(1.875,0.5)(2.5,0)}
\put(2,6){\VertexkTupleFib{1}{4}{7}{-0.7}{0.05}}
\put(3.25,5){\VertexkTupleFib{1}{4}{8}{0}{0.1}}
\put(4.5,4){\VertexkTupleFib{1}{4}{9}{0.1}{-0.05}}
\put(2,5){\VertexkTupleFib{1}{5}{7}{-0.7}{0.05}}
\put(3.25,4){\VertexkTupleFib{1}{5}{8}{0.05}{0.1}}
\put(4.5,3){\VertexkTupleFib{1}{5}{9}{0.1}{-0.05}}
\put(3.25,3){\VertexkTupleFib{1}{6}{8}{-0.725}{0.0}}
\put(4.5,2){\VertexkTupleFib{1}{6}{9}{0.1}{-0.05}}
\put(1,5){\VertexkTupleFib{2}{4}{7}{-0.7}{0.05}}
\put(2.25,4){\VertexkTupleFib{2}{4}{8}{-0.7}{-0.2}}
\put(3.5,3){\VertexkTupleFib{2}{4}{9}{0.07}{-0.1}}
\put(1,4){\VertexkTupleFib{2}{5}{7}{-0.7}{0.05}}
\put(2.25,3){\VertexkTupleFib{2}{5}{8}{-0.775}{-0.05}}
\put(3.5,2){\VertexkTupleFib{2}{5}{9}{0.075}{-0.1}}
\put(2.25,2){\VertexkTupleFib{2}{6}{8}{-0.7}{0}}
\put(3.5,1){\VertexkTupleFib{2}{6}{9}{0.1}{-0.1}}
\put(0,3){\VertexkTupleFib{3}{5}{7}{-0.75}{-0.05}}
\put(1.25,2){\VertexkTupleFib{3}{5}{8}{-0.7}{-0.15}}
\put(2.5,1){\VertexkTupleFib{3}{5}{9}{-0.7}{-0.15}}
\put(1.25,1){\VertexkTupleFib{3}{6}{8}{-0.7}{-0.05}}
\put(2.5,0){\VertexkTupleFib{3}{6}{9}{0.15}{-0.05}}
\put(1,5){\NEEdgeLabelForLatticeI{{\em 1}}}
\put(2,3.75){\NEEdgeLabelForLatticeI{{\em 1}}}
\put(2.1,2.85){\NEEdgeLabelForLatticeI{{\em 1}}}
\put(2.1,1.85){\NEEdgeLabelForLatticeI{{\em 1}}}
\put(3.75,3.25){\NEEdgeLabelForLatticeI{{\em 1}}}
\put(3.75,2.25){\NEEdgeLabelForLatticeI{{\em 1}}}
\put(3.6,1.1){\NEEdgeLabelForLatticeI{{\em 1}}}
\put(1.3,4.3){\NEEdgeLabelForLatticeI{{\em 1}}}
\put(2,5){\VerticalEdgeLabelForLatticeI{{\em 2}}}
\put(1,4){\VerticalEdgeLabelForLatticeI{{\em 2}}}
\put(3.25,4){\VerticalEdgeLabelForLatticeI{{\em 2}}}
\put(2.25,3){\VerticalEdgeLabelForLatticeI{{\em 2}}}
\put(4.5,3){\VerticalEdgeLabelForLatticeI{{\em 2}}}
\put(3.5,2){\VerticalEdgeLabelForLatticeI{{\em 2}}}
\put(3.25,3){\VerticalEdgeLabelForLatticeI{{\em 1}}}
\put(2.25,2){\VerticalEdgeLabelForLatticeI{{\em 1}}}
\put(1.25,1){\VerticalEdgeLabelForLatticeI{{\em 1}}}
\put(4.5,2){\VerticalEdgeLabelForLatticeI{{\em 1}}}
\put(3.5,1){\VerticalEdgeLabelForLatticeI{{\em 1}}}
\put(2.5,0){\VerticalEdgeLabelForLatticeI{{\em 1}}}
\put(3.15,5){\NWEdgeLabelForLatticeI{{\em 1}}}
\put(2.4,3.8){\NWEdgeLabelForLatticeI{{\em 1}}}
\put(2.9,4.2){\NWEdgeLabelForLatticeI{{\em 1}}}
\put(2.025,3.1){\NWEdgeLabelForLatticeI{{\em 1}}}
\put(1.025,2.1){\NWEdgeLabelForLatticeI{{\em 1}}}
\put(4.525,3.9){\NWEdgeLabelForLatticeI{{\em 2}}}
\put(4.15,3.2){\NWEdgeLabelForLatticeI{{\em 2}}}
\put(4.65,1.8){\NWEdgeLabelForLatticeI{{\em 2}}}
\put(3.525,1.9){\NWEdgeLabelForLatticeI{{\em 2}}}
\put(2.1,3.3){\NEEdgeLabelForLatticeI{{\em 2}}}
\put(2.1,1.3){\NEEdgeLabelForLatticeI{{\em 2}}}
\put(1.6,0.9){\NEEdgeLabelForLatticeI{{\em 2}}}
\put(1.225,0.2){\NEEdgeLabelForLatticeI{{\em 2}}}
\put(0.1,3.1){\NEEdgeLabelForLatticeI{{\em 2}}}
\put(1.25,2){\NEEdgeLabelForLatticeI{{\em 2}}}
\put(2.45,1.3){\NWEdgeLabelForLatticeI{{\em 2}}}
\put(2.3,0.8){\NEEdgeLabelForLatticeI{{\em 2}}}
\put(3.5,0.1){\NWEdgeLabelForLatticeI{{\em 2}}}
\end{picture}
\end{center}
\end{figure}

Versions of symmetric Fibonaccian lattices -- although not by this name -- appear in the 1976 paper \cite{BK} as Examples 2.1 and 2.2, but only for $n=3$.  
In the 1982 paper \cite{Gan}, Gansner considers versions of SFLs for general $n$. 
In particular, Gansner's lattices are (uncolored) distributive lattices of order ideals from `up-down' posets there denoted $Q(n-1,k(n-1))$. 
SFLs with $n=3$ are re-introduced in an explicitly Fibonaccian context in the 1985 paper \cite{HH}. 
Also see Exercise 3.23a of \cite{StanleyText} for SFLs with $n=3$.  
In the computing literature, SFLs for $n=3$ are called `Fibonacci cubes' and are considered as uncolored and undirected graphs, see \cite{Hsu}. 
See \cite{LotsaPs} for an interesting application of SFLs (in a different guise, and for $n=3$ only) to topology. 
The paper \cite{Munarini--Salvi} considers various aspects of the sequence of rank sizes for the $n=3$ version of SFLs and some related `Lucas' lattices; it would be interesting to see which of those results nicely extend to SFLs for general $n$. 
The above papers also consider distributive lattices with sizes $1, 2, 5, 13, 34, 89, \ldots$ comprising the other `half' of the Fibonacci sequence, but these lattices are not rank symmetric and are not considered here. 

A principal contribution of this paper is to realize each symmetric Fibonaccian lattice $L^{\mbox{\tiny Fib}}(n,k)$ as a model for a certain representation of the special linear Lie algebra $\mathfrak{sl}(n,\mathbb{C})$. 
We precede the development of these results with a brief study, in Section 2 below, of some enumerative aspects of SFLs, culminating in \EnumerativeTheorem\ and \EnumerativeCorollary.  
This work is aided by our introduction of some new recursively-defined triangular arrays that will be used in Section 5 to describe the rank sequences of SFLs.  
Many of these enumerative ideas are generalized in \cite{DDHK}. 
In Section 3, we provide a readable, first-principles pr\'{e}cis of key combinatorial notions related to $\mathfrak{sl}(n,\mathbb{C})$-representation theory.  
The relationship between the foregoing ideas and SFLs is developed in Sections 4 and 5. 
In particular, we show in \FibGTTheorem/\FibGTCorollary\ that $L^{\mbox{\tiny Fib}}(n,k)$ coincides with a certain skew-tabular lattice $L_{\mytinyA_{n-1}}^{\mbox{\tiny skew}}(\mysmallP/\mysmallQ)$ as in \cite{DD}. 
This skew-tabular lattice encodes, in a certain way, the actions of $\mathfrak{sl}(n,\mathbb{C})$-generators on a weight basis and has as its associated weight-generating function the skew Schur function we denote $\vartheta_{_{\mytinyP/\mytinyQ}}$. 
In \RGFTheorem, we produce two new explicit descriptions of the rank generating function $\RGF(L^{\mbox{\tiny Fib}}(n,k);q)$. 
One description employs the new recursively-defined triangular arrays of integers from Section 2. 
The other description borrows from \cite{DD} by realizing $\RGF(L^{\mbox{\tiny Fib}}(n,k);q)$ as a sum of rational expressions taken over a certain set $\mathcal{B}^{\mbox{\tiny Fib}}(n,k)$ of distinguished elements from $L^{\mbox{\tiny Fib}}(n,k)$. 
In Section 6, we consider some related problems which might be of interest. 
First, for all integers $n \geq 3$ and for all $k \in \{1,2,3\}$, we argue that each $L^{\mbox{\tiny Fib}}(n,k)$ enjoys a particular extremal property (the `edge-minimal' property) as a model for its associated $\mathfrak{sl}(n,\mathbb{C})$-representation and also appears to enjoy a certain uniqueness property (the `solitary' property); we believe $L^{\mbox{\tiny Fib}}(n,k)$ affords these same properties when $k > 3$. 
Second, we conjecture that the set $\mathcal{B}^{\mbox{\tiny Fib}}(3,k)$ is equinumerous with the set of what we call `topside peakless Motzkin paths' of length $k$ (cf.\ OEIS-A004148).

\vspace*{0.1in}
\noindent {\bf \S \ArraySection\ Some Fibonacci-related arrays and the enumeration of SFL cardinalities.}  In this section, we present some triangular arrays of nonnegative integers whose row sums are related to the Fibonacci numbers and to certain of their generalizations. 
Most importantly, we use the rows of these arrays to define univariate polynomials from which the cardinalities of the symmetric Fibonaccian lattices can be recovered in (perhaps) somewhat surprising ways. 

The first of the arrays we consider is well known and can be viewed as a `skewing' of Pascal's triangle to depict its Fibonacci diagonals as rows:  
{
\begin{center}
\renewcommand{\baselinestretch}{1.2}
\footnotesize
\begin{tabular}{ccccc}
1 & & & & \\
1 & & & & \\
1 & 1 & & & \\
1 & 2 & & & \\
1 & 3 & 1 & & \\
1 & 4 & 3 & & \\
1 & 5 & 6 & 1 & \\
1 & 6 & 10 & 4 & \\
1 & 7 & 15 & 10 & 1 \\
1 & 8 & 21 & 20 & 5 
\end{tabular}

\vspace*{0.1in}
ETC.
\end{center}
}
\noindent 
This {\em Fibonacci right-triangular array} is herein denoted as $\mymathfrak{F}$, and the $j^{\mbox{\tiny \em th}}$ entry of the $k^{\mbox{\tiny \em th}}$ row is denoted $f_{k,j}$, so that \[\mymathfrak{F} = \left(f_{k,j}\, \mbox{\small $= {k-j \choose j}$}\right)_{k\in\{0,1,2,\ldots\}, j\in\{0,1,\ldots,\lfloor\frac{k}{2}\rfloor\}}.\] 
We can build $\mymathfrak{F}$ recursively by declaring $f_{0,0} :=  1 =: f_{1,0}$ and then applying the recurrence 
\[f_{k,j} := f_{k-1,j} + f_{k-2,j-1}\] 
for integers $j$ and $k$ with $k \geq 2$ and with the understanding that $f_{k,j} := 0$ when $j < 0$ or $j > \lfloor \frac{k}{2} \rfloor$. 
Using this defining recurrence, one can see that for any nonnegative integer $k$ we have $\displaystyle \sum_{j=0}^{\lfloor \frac{k}{2} \rfloor}f_{k,j} = f_{k}$.  
The $k^{\mbox{\tiny \em th}}$ row entries of $\mymathfrak{F}$ have many combinatorial interpretations (e.g.\ see the comments / references of OEIS-A011973) and can be thought of as a refinement of the $k^{\mbox{\tiny \em th}}$ Fibonacci number. 

Our main interest in $\mymathfrak{F}$ is in a family of polynomials whose (unsigned) coefficients come from this array and which we call {\em sign-alternating Fibonacci array polynomials}: 
\[\altFibpoly{k}(x) := \sum_{j=0}^{\lfloor \frac{k}{2} \rfloor} (-1)^{j}f_{k,j}\, x^{\lfloor \frac{k}{2} \rfloor - j},\] 
with $\altFibpoly{-1}(x) := 0$ when needed. 
(The use of the overline in the notation `$\altFibpoly{k}(x)$' serves as a reminder that the latter polynomial is at least partly distinguished from many other so-called Fibonacci polynomials by the presence of alternating signs.) 
From the defining recurrence of $\mymathfrak{F}$ we get the following recurrence relation for all integers $k \geq 2$:  
\[\altFibpoly{k}(x) = x^{(k-1)\, \mbox{\scriptsize mod 2}}\,\altFibpoly{k-1}(x) - \altFibpoly{k-2}(x),\]
taking $\altFibpoly{0}(x) = 1 = \altFibpoly{1}(x)$. 
In \EnumerativeTheorem, we connect these polynomials to our SFLs. 

Next, we consider a family of arrays indexed by integers $n \geq 2$. 
For $n > 3$, we believe these integer arrays are new. 
The $n^{\mbox{\tiny \em th}}$ {\em symmetric Fibonacci triangle} $\mymathfrak{A}^{(n)} = (\mya^{(n)}_{k,r})$ is defined recursively as follows. 
For each nonnegative integer $k$, set $\mathcal{R}_{n,k} := \{-k(n-1),-k(n-1)+2,\ldots,k(n-1)-2,k(n-1)\}$, and for $r \not\in \mathcal{R}_{n,k}$ set $\mya^{(n)}_{k,r} := 0$. 
For initial values, take $\mya^{(n)}_{0,0} := 1$ with $\mya^{(n)}_{1,r} := 1$ for each $r \in \mathcal{R}_{n,1}$. 
Then for $k \geq 2$ and $j \in \mathbb{Z}$, let $\displaystyle \mya^{(n)}_{k,r} := \left(\sum_{s \in \mathcal{R}_{n,1}}\mya^{(n)}_{k-1,r+s}\right) - \mya^{(n)}_{k-2,r}$. 
Here, for example, is part of $\mymathfrak{A}^{(4)}$:  
{
\begin{center}
\renewcommand{\baselinestretch}{1.2}
\footnotesize
\begin{tabular}{ccccccccccccccccccc}
 & & & & & & & & & 1 & & & & & & & & & \\
 & & & & & & 1 & & 1 & & 1 & & 1 & & & & & & \\
 & & & 1 & & 2 & & 3 & & 3 & & 3 & & 2 & & 1 & & & \\
1 & & 3 & & 6 & & 8 & & 10 & & 10 & & 8 & & 6 & & 3 & & 1 
\end{tabular}

\vspace*{0.1in}
ETC.
\end{center}
}
\noindent
We declare that $\displaystyle {A}^{(n)}_{k}(x) := \sum_{r \in \mathcal{R}_{n,k}}\mya^{(n)}_{k,r}x^{\frac{1}{2}\left(\rule[-1.25mm]{0mm}{3.5mm}k(n-1)-r\right)}$, and set $\displaystyle {A}^{(n)}_{-1}(x) := 0$. 

For $n=3$, the symmetric Fibonacci triangle $\mymathfrak{A}^{(3)}$, which we found in \cite{KK}, is known mainly as the subarray of even numbered rows of the triangle appearing in the OEIS as entry A079487. 
In \cite{KK}, it is shown that the rows of the OEIS-A079487 array are the rank numbers of the lattice of order ideals of zigzag posets. 
Thus, the $k^{\mbox{\tiny th}}$ even numbered row of the OEIS-A079487 array is the list of rank numbers of our symmetric Fibonaccian lattice $L^{\mbox{\tiny Fib}}(3,k)$. 
This latter fact is re-proven here as a special case of \RGFTheorem. 

The main result of this section, \EnumerativeTheorem\ below, is enumerative. 
Part of its proof involves a careful analysis of the ranks of $L^{\mbox{\tiny Fib}}(n,k)$.  
This lattice has maximal element $M = (1,n+1,\ldots,(k-1)n+1)$, minimal element $N = (n,2n,\ldots,kn)$, and length $k(n-1)$. 
Then, one can see that the rank function $\rho: L^{\mbox{\tiny Fib}}(n,k) \longrightarrow \{0,\ldots,k(n-1)\}$ is given by \[\rho(T) = k(n-1) - \sum_{i=1}^{k}(T_{i}-(i-1)n-1) = \frac{1}{2}k(k+1)n - \sum_{i=1}^{k} T_{i}.\] 
Define the {\em centered rank function} $\rho^{\mbox{\tiny ctr}}: L^{\mbox{\tiny Fib}}(n,k) \longrightarrow \mathcal{R}_{n,k}$ by the rule \[\rho^{\mbox{\tiny ctr}}(T) = 2\rho(T)-k(n-1) = k(kn+1) - 2\sum_{i=1}^{k} T_{i}.\] 
We set $\mathcal{C}_{r}(n,k) := (\rho^{\mbox{\tiny ctr}})^{-1}(r)$ when $r \in \mathcal{R}_{n,k}$; otherwise set $\mathcal{C}_{r}(n,k) := \emptyset$. 

\noindent
{\bf \RankProposition}\ \ {\sl Let $n$ and $k$ be positive integers, with $n \geq 2$.  For all integers $r$,}
\[\left|\rule[-1.5mm]{0mm}{4.75mm}\mathcal{C}_{r}(n,k)\, \mbox{\large$\sqcup$}\, \mathcal{C}_{r}(n,k-2)\right| = \left|\rule[-1.5mm]{0mm}{4.75mm}\bigsqcup_{s \in \mathcal{R}_{n,1}} \mathcal{C}_{r+s}(n,k-1)\right|,\]
{\sl where `$\, \sqcup$' denotes disjoint union.  
Moreover, $\displaystyle \left|\rule[-1.5mm]{0mm}{4.75mm}\mathcal{C}_{r}(n,k)\right| = \mya^{(n)}_{k,r}$.}

{\em Proof.} Our proof is bijective. 
Let $S \in \bigsqcup_{s \in \mathcal{R}_{n,1}} \mathcal{C}_{r+s}(n,k-1)$, so $r+s = \rho^{\mbox{\tiny ctr}}(S) = (k-1)[(k-1)n+1-2\sum_{i=1}^{k-1}S_{i}$ for some $s \in \mathcal{R}_{n,1}$. 
It is clear that there is exactly one possible number $T'_{k}$ such that, when we take $T' := (T'_{1},\ldots,T'_{k-1},T'_{k})$ with $T'_{i} = S_{i}$ for $1 \leq i \leq k-1$, then $r = k(kn+1) - 2\sum_{i=1}^{k} T'_{i}$, namely $T'_{k} = \frac{1}{2}(s+2nk-n+1)$. 
Now, by routine manipulation of inequalities, one can see that $T'_{k} \in \{(k-1)n+1,(k-1)n+2,\ldots,kn\}$ if and only if $s \in \mathcal{R}_{n,1}$, and of course the latter is true. 
Moreover, $T'_{k} = (k-1)n+1$ if and only if $s = -(n-1)$. 
We now consider two cases: (1) $T_{k-1} < (k-1)n$ or $s > -(n-1)$ and (2) $T_{k-1} = (k-1)n$ and $s = -(n-1)$. 
In case (1), the $k$-tuple $T'$ is an element of $\mathcal{C}_{r}(n,k)$, and we rename $T'$ as simply $T$. 
In case (2), the $k$-tuple $T'$ fails the requirement that $T'_{k-1}+1 \ne T'_{k}$, so we let $T := (S_{1},\ldots,S_{k-2})$, which is in $L^{\mbox{\tiny Fib}}(n,k-2)$. 
A routine calculation shows that $\rho^{\mbox{\tiny ctr}}(T) = r$, hence $T \in \mathcal{C}_{r}(n,k-2)$. 
We now have a set mapping $\bigsqcup_{s \in \mathcal{R}_{n,1}} \mathcal{C}_{r+s}(n,k-1) \longrightarrow \mathcal{C}_{r}(n,k)\, \mbox{\large$\sqcup$}\, \mathcal{C}_{r}(n,k-2)$ given by $S \mapsto T$. 
It is clear that our procedure for producing $T$ can be reversed in order to recover $S$, so our set mapping is a bijection.

To see that the numbers $\displaystyle \left|\rule[-1.5mm]{0mm}{4.75mm}\mathcal{C}_{r}(n,k)\right|$ satisfy the same recurrence as elements of the array $\mymathfrak{A}^{(n)} = (\mya^{(n)}_{k,r})$, 
rewrite the just-established identity $\left|\rule[-1.5mm]{0mm}{4.75mm}\mathcal{C}_{r}(n,k)\, \mbox{\large$\sqcup$}\, \mathcal{C}_{r}(n,k-2)\right| = \left|\rule[-1.5mm]{0mm}{4.75mm}\bigsqcup_{s \in \mathcal{R}_{n,1}} \mathcal{C}_{r+s}(n,k-1)\right|$ 
as $\left|\rule[-1.5mm]{0mm}{4.75mm}\mathcal{C}_{r}(n,k)\right| = \left|\rule[-1.5mm]{0mm}{4.75mm}\bigsqcup_{s \in \mathcal{R}_{n,1}} \mathcal{C}_{r+s}(n,k-1)\right| - \left|\rule[-1.5mm]{0mm}{4.75mm}\mathcal{C}_{r}(n,k-2)\right|$. 
We conclude that $\displaystyle \left|\rule[-1.5mm]{0mm}{4.75mm}\mathcal{C}_{r}(n,k)\right| = \mya^{(n)}_{k,r}$.\hfill\QED

\noindent 
{\bf \EnumerativeTheorem}\ \ {\sl Let $n$ and $k$ be positive integers, with $n \geq 2$.  Set} 
$\mathscr{L}_{n,k} := \left|\rule[-1.5mm]{0mm}{4.75mm}L^{\mbox{\tiny Fib}}(n,k)\right|$, 
$\mathscr{F}_{n,k} := n^{k\, \mbox{\scriptsize mod 2}}\, \overline{F}_{k}(n^{2})$, 
{\sl and} $\mathscr{A}_{n,k} := {A}^{(n)}_{k}(1)$. 
{\sl Then for each $\mathscr{X} \in \{\mathscr{L},\mathscr{F},\mathscr{A}\}$ we have $\mathscr{X}_{n,k} = n\mathscr{X}_{n,k-1} - \mathscr{X}_{n,k-2}$, with $\mathscr{X}_{n,-1} = 0$ and $\mathscr{X}_{n,0} = 1$. 
In particular, $\mathscr{L}_{n,k} = \mathscr{F}_{n,k} = \mathscr{A}_{n,k}$.}

{\em Proof.} It is routine to verify that the $\mathscr{F}_{n,k}$'s satisfy the claimed recurrence. 
That $\mathscr{A}_{n,k} = n\mathscr{A}_{n,k-1} - \mathscr{A}_{n,k-2}$ is easily established using the defining recurrence for the $\mya^{(n)}_{k,j}$'s. 
By \RankProposition, $\displaystyle \left|\rule[-1.5mm]{0mm}{4.75mm}\mathcal{C}_{r}(n,k)\right| = \mya^{(n)}_{k,r}$, from which it follows that $\mathscr{L}_{n,k} := \left|\rule[-1.5mm]{0mm}{4.75mm}L^{\mbox{\tiny Fib}}(n,k)\right| = \sum_{r \in \mathcal{R}_{n,k}}\left|\rule[-1.5mm]{0mm}{4.75mm}\mathcal{C}_{r}(n,k)\right| = \sum_{r \in \mathcal{R}_{n,k}}\mya^{(n)}_{k,r} = {A}^{(n)}_{k}(1) = \mathscr{A}_{n,k}$.\hfill\QED

For a fixed integer $n \geq 2$, \EnumerativeTheorem\ offers several interpretations and contexts for the sequence of numbers $\{\mathscr{X}_{n,k}\}_{k \geq -1}$ defined by the recurrence $\mathscr{X}_{n,k} = n\mathscr{X}_{n,k-1} - \mathscr{X}_{n,k-2}$ when $k \geq 1$, with $\mathscr{X}_{n,-1} = 0$ and $\mathscr{X}_{n,0} = 1$. 
The number $\mathscr{X}_{n,k} = \mathscr{L}_{n,k}$ enumerates the elements of the SFL $L^{\mbox{\tiny Fib}}(n,k)$. 
A consequence of \RankProposition\ is a result that we formally record in \RGFTheorem, namely that the polynomial ${A}^{(n)}_{k}(q)$ is the rank generating function of $L^{\mbox{\tiny Fib}}(n,k)$ and therefore, by \EnumerativeTheorem, serves as a $q$-analog of the number $\mathscr{X}_{n,k} = \mathscr{A}_{n,k}$. 
Another connection, observed by the anonymous referee, is that the defining recurrence for the sequence $\{\mathscr{X}_{n,k}\}_{k \geq -1}$, and another closely related recursively defined sequence, was introduced in \cite{KM} to help analyze the real and especially imaginary roots of the rank two symmetric hyperbolic Kac--Moody Lie algebra associated with the generalized Cartan matrix {\footnotesize $\left(\!\!\!\begin{array}{cc}2 & -n\\ -n & 2\end{array}\!\!\!\right)$}. 
These recurrences are used again in \cite{KLL} to bound imaginary root multiplicities. 
An interesting question is whether/how such studies of root systems for rank two symmetric hyperbolic Kac--Moody Lie algebras might connect with the representation theoretic context of this paper.   

We close this section with a brief investigation of some connections between \EnumerativeTheorem\ and some integer sequences from the OEIS.  
Notice that in the statement of \EnumerativeTheorem, the recurrence $\mathscr{X}_{n,k} = n\mathscr{X}_{n,k-1} - \mathscr{X}_{n,k-2}$, with $\mathscr{X}_{n,-1} = 0$ and $\mathscr{X}_{n,0} = 1$ is defined for $n=1$, yielding the period six repeating sequence $(\mathscr{X}_{1,k})_{k\geq 0} = (1, 1, 0, -1, -1, 0, 1, 1, 0, -1, -1, 0, \ldots)$. 
The definitions of $\mathscr{F}_{n,k}$ and $\mathscr{A}_{n,k}$ make sense when $n=1$ and agree with the preceding sequence.
A table of small values for the $\mathscr{F}_{n,k}$'s/$\mathscr{A}_{n,k}$'s is given in \TableFig. 
\begin{figure}[ht]
\begin{center}
\begin{tabular}{|c||c|c|c|c|c|c|c|c|}
\hline
$n$ {\tt \char`\\ } $k$ & 0 & 1 & 2 & 3 & 4 & 5 & 6 & $\cdots$\\
\hline
\hline
1 & 1 & 1 & 0 & -1 & -1 & 0 & 1 & $\cdots$\\
\hline
2 & 1 & 2 & 3 & 4 & 5 & 6 & 7 & $\cdots$\\
\hline
3 & 1 & 3 & 8 & 21 & 55 & 144 & 377 & $\cdots$\\
\hline
4 & 1 & 4 & 15 & 56 & 209 & 780 & 2911 & $\cdots$\\
\hline
5 & 1 & 5 & 24 & 115 & 551 & 2640 & 12,649 & $\cdots$\\
\hline
6 & 1 & 6 & 35 & 204 & 1189 & 6930 & 40,391 & $\cdots$\\
\hline
7 & 1 & 7 & 48 & 329 & 2255 & 15,456 & 105,937 & $\cdots$\\
\hline
8 & 1 & 8 & 63 & 496 & 3905 & 30,744  & 242,047 & $\cdots$\\
\hline
9 & 1 & 9 & 80 & 711 & 6319  & 56,160 & 499,121 & $\cdots$\\
\hline
10 & 1 & 10 & 99 & 980 & 9701 & 96,030 & 950,599 & $\cdots$\\
\hline
$\vdots$ & $\vdots$ & $\vdots$ & $\vdots$ & $\vdots$ & $\vdots$ & $\vdots$ & $\vdots$ & $\ddots$\\
\hline
\end{tabular}

\vspace*{0.15in}
{\bf \TableFig:} A table of small values for the $\mathscr{F}_{n,k}$'s/$\mathscr{A}_{n,k}$'s.
\end{center}
\vspace*{-0.25in}
\end{figure}
In this array, the sequence of antidiagonal entries begins with {\small $1, 1, 1, 1, 2, 0, 1, 3, 3, -1, 1, 4, 8, 4, -1, \ldots$} and appears to be the same as the rectangular array of OEIS  sequence A179943 but modified to include OEIS-A010892 as a first row. 
As a consequence of the preceding theorem, we can formally connect our work here to OEIS-A179943 and OEIS-A179944. 
 
\noindent
{\bf \EnumerativeCorollary}\ \ {\sl Let $n$ and $k$ be integers with $n \geq 1$ and $k \geq -1$. 
Following OEIS-A179943, let $\mathscr{M}_{n,k}$ be the $(2,1)$ entry of the $2 \times 2$ matrix} {\footnotesize $\left(\begin{array}{cc}1 & n-2\\ 1 & n-1\end{array}\right)^{k+1}$.} {\sl Then $\mathscr{M}_{n,k} = \mathscr{F}_{n,k}=\mathscr{A}_{n,k}$.
Moreover, in the above rectangular array of $\mathscr{M}_{n,k}$'s/$\mathscr{F}_{n,k}$'s/$\mathscr{A}_{n,k}$'s, the sum of the entries on the $d^{\mbox{\tiny th}}$ antidiagonal (beginning with $d=0$) is}
\[\sum_{c=0}^{d} (c+1)^{(d-c)\, \mbox{\scriptsize mod}\, 2}\overline{F}_{d-c}((c+1)^{2}),\]
{\sl which coincides with the sum of the $d^{\mbox{\tiny th}}$ term of OEIS-A010892 and the $(d-1)^{\mbox{\tiny st}}$ term of OEIS-A179944. (In the latter sequence, we regard the $(-1)^{\mbox{\tiny st}}$ term to be zero.)}

{\em Proof.} Fix $n$. 
Clearly $\mathscr{M}_{n,-1} = 0$ and $\mathscr{M}_{n,0} = 1$. 
A simple induction argument shows that for $k \geq 1$, we get $\mathscr{M}_{n,k} = n\mathscr{M}_{n,k-1} - \mathscr{M}_{n,k-2}$. 
The displayed expression in the concluding statement of the corollary just interprets entries of the $d^{\mbox{\tiny th}}$ antidiagonal as certain $\mathscr{F}_{n,k}$'s. 
The claims relating this antidiagonal sum to the terms of the OEIS sequences A010892 and A179944 follow from the fact that the $\mathscr{M}_{n,k}$'s satisfy the recurrence noted at the beginning of this proof.\hfill\QED

\vspace*{0.1in}
\noindent {\bf \S \SetupSection\ Some first-principles background on relevant combinatorial aspects of $\mathfrak{sl}(n,\mathbb{C})$ representation theory.} 
For notation and language regarding the relevant order- and representation-theoretic concepts, we largely follow \cite{DonDiamond} and \cite{DD}. 
For more general background on Lie algebras/groups and their representations, \cite{Hum} and \cite{FH} are two standard references. 
Here, we offer a distilled version of these supporting ideas in order make more accessible the main results of the following sections (\TheoremList). 
This section can be lightly browsed on a first reading. 

Continue with the assumption that $n$ is an integer greater than one.  
It is helpful to define $E_{i,j}$ to be an $n \times n$ matrix with a `$1$' in the $(i,j)$ position and $0$'s elsewhere. 
The {\em special linear Lie algebra} $\mathfrak{sl}(n,\mathbb{C})$ is the complex vector space consisting of the traceless $n \times n$ complex matrices together with the binary Lie bracket operation defined by $[A,B] := AB-BA$ that is {\em Jacobi-associative} in that $[A,[B,C]] + [B,[C,A]] + [C,[A,B]] = 0$ when $A, B, C \in \mathfrak{sl}(n,\mathbb{C})$.  
This Lie algebra is {\em simple} in that it contains no nonzero and proper {\em ideal} $\mathfrak{I}$, where the latter is a subspace that satisfies $[A,B] \in \mathfrak{I}$ when $A \in \mathfrak{sl}(n,\mathbb{C})$ and $B \in \mathfrak{I}$.  
Moreover, $\mathfrak{sl}(n,\mathbb{C})$ possesses a distinguished set of generators $\{X_{i},Y_{i},H_{i}\}_{i=1}^{n-1}$ with $X_{i} := E_{i,i+1}$, $Y_{i} := E_{i+1,i}$, and $H_{i} := E_{i,i}-E_{i+1,i+1}$. 
For fixed $i$, each triple $(X_{i},Y_{i},H_{i})$ generates a Lie subalgebra isomorphic to $\mathfrak{sl}(2,\mathbb{C})$. 
Note that each $H_{i}$ is diagonal.

Within the classification of finite-dimensional simple Lie algebras over $\mathbb{C}$, our special linear Lie algebra $\mathfrak{sl}(n,\mathbb{C})$ has type `$\myA_{n-1}$'. 
Data that allows us to construct, via generators and relations, an isomorphic version of $\mathfrak{sl}(n,\mathbb{C})$ is contained in the associated $(n-1) \times (n-1)$ {\em Cartan matrix}
\[(a_{i,j})_{i,j \in \{1,2,\ldots,n-1\}} := \left(\begin{array}{ccccccc}
2 & -1 & 0 & \cdots & 0 & 0 & 0\\
-1 & 2 & -1 & \cdots & 0 & 0 & 0\\
0 & -1 & 2 & \cdots & 0 & 0 & 0\\
\vdots & \vdots & \vdots & \ddots & \vdots & \vdots & \vdots\\
0 & 0 & 0 & \cdots & 2 & -1 & 0\\
0 & 0 & 0 & \cdots & -1 & 2 & -1\\
0 & 0 & 0 & \cdots & 0 & -1 & 2
\end{array}\right).\] 
Let $\Omega := \{\omega_{i}\}_{i=1}^{n-1}$ be the standard basis for the $\mathbb{Z}$-module $\mathbb{Z}^{n-1}$, and let $\alpha_{i}$ denote the $i^{\mbox{\tiny th}}$ row vector of our Cartan matrix, so $\alpha_{i} := \sum_{j=1}^{n-1}a_{i,j}\omega_{j}$.  
Using this matrix, we define a Lie algebra over $\mathbb{C}$ with the following presentation by generators and relations:  
\[\mathfrak{g}(\myA_{n-1}) := \langle \myqx_{i},\myqy_{i},\myqh_{i}\, |\, \mbox{relations (1), (2), (3), (4), (5) below}\rangle,\]
where for all $i,j \in \{1,2,\ldots,n-1\}$ we have (1) $[\myqx_{i},\myqy_{j}] = \delta_{i,j}\myqh_{i}$, (2) $[\myqh_{i},\myqx_{j}] = a_{i,j}\myqx_{j}$, (3) $[\myqh_{i},\myqy_{j}] = -a_{i,j}\myqy_{j}$, (4) $\mbox{ad}(\myqx_{i})^{1-a_{i,j}}(\myqx_{j}) = 0$ when $i \ne j$, and (5) $\mbox{ad}(\myqy_{i})^{1-a_{i,j}}(\myqy_{j}) = 0$ when $i \ne j$, with $\mbox{ad}(\myqz)(\myqz')$ understood to be $[\myqz,\myqz']$ for any $\myqz,\myqz'$ in $\mathfrak{g}(\myA_{n-1})$. 
It can be seen that  $\mathfrak{g}(\myA_{n-1})$ is isomorphic to $\mathfrak{sl}(n,\mathbb{C})$ via the Lie algebra homomorphism induced by $\myqx_{i} \mapsto X_{i}$, $\myqy_{i} \mapsto Y_{i}$, and $\myqh_{i} \mapsto H_{i}$. 
For $i \in \{1,\ldots,n-1\}$, let $S_{i}: \mathbb{Z}^{n-1} \longrightarrow \mathbb{Z}^{n-1}$ be the order two invertible linear transformation given by $S_{i}(\mu) = \mu -\mu_{i}\alpha_{i}$ when $\mu = \sum_{j =1}^{n-1}\mu_{j}\omega_{j}$. 
Identify each $S_{i}$ with its matrix representation in the basis $\Omega$, so $\{S_{i}\}_{i = 1}^{n-1} \subseteq GL_{n-1}(\mathbb{Z})$, where the latter is the group of $(n-1) \times (n-1)$ unimodular matrices.  
The {\em Weyl group} $W(\myA_{n-1})$ associated with $\mathfrak{g}(\myA_{n-1}) \cong \mathfrak{sl}(n,\mathbb{C})$ is the subgroup of $GL_{n-1}(\mathbb{Z})$ generated by $\{S_{i}\}_{i = 1}^{n-1}$. 

A finite-dimensional complex vector space $V$ is a $\mathfrak{g}(\myA_{n-1})${\em -module} when there is an action `$\myqz.v$' of special linear Lie algebra elements $\myqz$ on vectors $v$ so that $[\myqz,\myqz'].v = \myqz.(\myqz'.v)-\myqz'.(\myqz.v)$ for all $\myqz, \myqz' \in \mathfrak{g}(\myA_{n-1})$ and all $v \in V$ and such that the action respects the defining relations (1) through (5) for $\mathfrak{g}(\myA_{n-1})$. 
We also call $V$ a {\em representing space} for, or a {\em representation} of, $\mathfrak{g}(\myA_{n-1})$. 
Given such a $\mathfrak{g}(\myA_{n-1})$-module $V$, it is always possible to find a basis, called a {\em weight basis}, such that the representing matrices for the $\myqh_{i}$'s are diagonal with integer entries. 
Suppose a weight basis $\{v_{\telt}\}_{\telt \in R}$ is indexed by some set $R$ of, say, combinatorial objects. 
Define $\mym_{i}(\telt) \in \mathbb{Z}$ by $\myqh_{i}.v_{\telt} = \mym_{i}(\telt)v_{\telt}$. 
The {\em weight} of $v_{\telt}$, or of the index $\telt$, is the vector of eigenvalues $\wt(v_{\telt}) = \wt(\telt) := \sum \mym_{i}(\telt)\omega_{i}$. 
Such a weight basis can be depicted graphically as follows. 
Regard the index set $R$ to be a set of vertices for an edge-colored directed graph where, for any $\selt$ and $\telt$ in $R$, we write $\selt \myarrow{i} \telt$ if, when we write $\myqx_{i}.v_{\selt} = \sum_{\uelt \in R}\myqX_{\uelt,\selt}v_{\uelt}$ and $\myqy_{i}.{\telt} = \sum_{\relt \in R}\myqY_{\relt,\telt}v_{\relt}$, at least one of $\myqX_{\telt,\selt}$ or $\myqY_{\selt,\telt}$ is nonzero. 
The edge-colored directed graph $R$ is the {\em supporting graph} for the given weight basis; the {\em representation diagram} for the weight basis is the edge-colored digraph $R$ together with the scalar pairs $(\myqX_{\telt,\selt},\myqY_{\selt,\telt})$ associated with the edges $\selt \myarrow{i} \telt$.  
On a representation diagram edge $\selt \myarrow{i} \telt$, the associated {\em edge product} is $\myqP_{\selt,\telt} := \myqX_{\telt,\selt}\myqY_{\selt,\telt}$. 

The supporting graph and representation diagram for a given weight basis have much combinatorial structure, some of which we note here. 
Let $\mu = \sum \mu_{i}\omega_{i}$ be a vector in $\mathbb{Z}^{n-1}$, and let $V_{\mu}$ be the subspace spanned by the set $\{v_{\telt}\, |\, \wt(\telt)=\mu\}$. 
Observe that if $v \in V_{\mu}$, then $\myqx_{i}.v \in V_{\mu+\alpha_{i}}$ and $\myqy_{i}.v \in V_{\mu-\alpha_{i}}$. 
So, if $\selt \myarrow{i} \telt$ in $R$, then $\wt(\selt)+\alpha_{i}=\wt(\telt)$. 
It can be deduced that $R$ is an acyclic directed graph that coincides with its own transitive reduction; that is, the {\em reachability relation} ($\selt \leq \telt$ in $R$ if and only if there is a directed path from $\selt$ to $\telt$) is a partial order on the elements of $R$ and the covering relations are exactly the directed edges of $R$. 
It can be further deduced that the poset $R$ is ranked by some surjective {\em rank function} $\rho: R \longrightarrow \{0,1,\ldots,l\}$ such that $\rho(\selt)+1=\rho(\telt)$ whenever $\selt \rightarrow \telt$ in $R$; in this case, the number $l$ is the length of $R$ with respect to $\rho$.  
The rank function is unique if $R$ is connected. 
We define the {\em depth function} $\delta: R \longrightarrow \{0,1,\ldots,l\}$ by $\delta(\telt) := l-\rho(\telt)$. 
The $i$-component $\comp_{i}(\telt)$ of some $\telt$ in $R$ is the set of all vertices in $R$ that can be reached from $\telt$ via some path of edges of color $i$ (where we disregard the direction of each edge) together with the set of directed edges from amongst all such possible paths. 
Then $\comp_{i}(\telt)$ is ranked with unique rank and depth functions denoted $\rho_{i}$ and $\delta_{i}$ respectively. 
One can easily check now that $\mym_{i}(\telt) = \rho_{i}(\telt)-\delta_{i}(\telt)$. 

Let $\{q\}$ and $\myvarZ := \{z_{1},\ldots,z_{n-1}\}$ be sets of indeterminants. 
The rank-generating function of $R$ is $\RGF(R;q) := \sum_{\telt \in R} q^{\rho(\telt)}$ and its weight-generating function is $\WGF(R;\myvarZ) := \sum_{\telt \in R}\myvarZ^{\smallwt(\telt)}$, where $\myvarZ^{\smallwt(\telt)} := z_{1}^{\mysmallerm_{1}(\telt)} \cdots z_{n-1}^{\mysmallerm_{n-1}(\telt)}$. 
The Weyl group $W(\myA_{n-1})$ acts on the $\mathbb{Z}$-module generated by the symbols $\{\myvarZ^{\mu}\}_{\mu \in \mathbb{Z}^{n-1}}$ by declaring $S_{i}.\myvarZ^{\mu} := \myvarZ^{S_{i}(\mu)}$ and extending $\mathbb{Z}$-linearly. 
It can be seen that for any supporting graph $R$ and any $i \in \{1,2,\ldots,n-1\}$, we have $S_{i}.\WGF(R;\myvarZ) = \WGF(R;\myvarZ)$. 
We say $\WGF(R;\myvarZ)$ is $W(\myA_{n-1})$-{\em invariant}. 
A consequence of the theory of so-called `characters' of Lie algebra representations is that, when $Q$ and $R$ are supporting graphs of weight bases for $\mathfrak{g}(\myA_{n-1})$-modules $V$ and $W$ respectively, then $V \cong W$ (in the natural sense of a $\mathfrak{g}(\myA_{n-1})$-module isomorphism) if and only if $\WGF(Q;\myvarZ) = \WGF(R;\myvarZ)$. 
That is, the $W(\myA_{n-1})$-invariant weight-generating function of a supporting graph uniquely identifies the associated representation. 

When a finite-dimensional $\mathfrak{g}(\myA_{n-1})$-module $V$ has no nontrivial and proper subspaces that are stable under the action, then $V$ is said to be {\em irreducible}. 
It follows from, say, the theory of Verma modules (see \cite{Hum}) that, when $V$ is irreducible, there exists a nonzero vector $v \in V$, unique up to scalar multiples, such that $\myqx_{i}.v = 0$ for all $i \in \{1,2,\ldots,n-1\}$. 
In this case, $v$ is necessarily a weight basis vector and its weight vector $\wt(v)$ is a nonnegative-integer linear combination of the vectors in $\Omega$. 
From here on, an integer linear combination $\sum a_{i}\omega_{i}$ with each $a_{i} \geq 0$ is a {\em dominant weight}. 
If $W$ is another finite-dimensional irreducible $\mathfrak{g}(\myA_{n-1})$-module whose unique maximal vector has weight coinciding with $\wt(v)$, then $W \cong V$. 
Moreover, it is known that for any dominant weight $\lambda = \sum a_{i}\omega_{i}$ of the vectors in $\Omega$, there exists an irreducible $\mathfrak{g}(\myA_{n-1})$-module whose unique maximal vector has weight $\lambda$. 
Within the theory of Verma modules, such a $\mathfrak{g}(\myA_{n-1})$-module `$V(\lambda)$' can be obtained via a certain quotient construction. 
{\em Complete reducibility} of semisimple Lie algebra representations is a major feature from the general theory that applies as follows here: Any finite-dimensional $\mathfrak{g}(\myA_{n-1})$-module $V$ is isomorphic to a direct sum $V(\lambda_{1}) \oplus \cdots \oplus V(\lambda_{l})$ of irreducible $\mathfrak{g}(\myA_{n-1})$-modules $V(\lambda_{k})$, for some list of dominant weights $\lambda_{1},\ldots,\lambda_{l}$. 
Moreover, if $R$ is any supporting graph for $V$ and if $R(\lambda_{k})$ is a generic supporting graph for $V(\lambda_{k})$, then $\WGF(R;\myvarZ) = \sum_{1 \leq k \leq l} \WGF(R(\lambda_{k});\myvarZ)$. 
Further, when $W \cong V(\mu_{1}) \oplus \cdots \oplus V(\mu_{m})$ for some dominant weights $\mu_{1},\ldots,\mu_{m}$, then $V \cong W$ if and only if $m=l$ and there is a permutation of the dominant weights $\mu_{1},\ldots,\mu_{m}$ such that each $\mu_{k}$ is just $\lambda_{k}$. 

Sometimes it is possible to start with a well-chosen set of combinatorial objects and then build from it a representation diagram without {\em a priori} knowledge of the associated weight basis. 
Here we note a particular such circumstance that is relevant for our work with SFL's. 
Say $L$ is a diamond-colored distributive lattice with edges colored by the set $\{1,\ldots,n-1\}$; let $\selt, \telt \in L$.  
The $i$-component $\comp_{i}(\telt)$ in $L$ has a unique rank (respectively, depth) function $\rho_{i}$ (resp.\ $\delta_{i}$).  
Without yet knowing whether $L$ is a supporting graph, the following quantities may still be defined as above: $\mym_{i}(\telt) := \rho_{i}(\telt)-\delta_{i}(\telt)$, $\wt(\telt) := \sum_{i=1}^{n-1}\mym_{i}(\telt)\omega_{i}$, $\myvarZ^{\smallwt(\telt)} := z_{1}^{\mysmallerm_{1}(\telt)} \cdots z_{n-1}^{\mysmallerm_{n-1}(\telt)}$, and $\WGF(L;\myvarZ) := \sum_{\telt \in L}\myvarZ^{\smallwt(\telt)}$. 
We say $L$ is $\myA_{n-1}$-{\em structured} if we have $\wt(\selt)+\alpha_{i} = \wt(\telt)$ whenever $\selt \myarrow{i} \telt$. 
Suppose now that we assign to each edge $\selt \myarrow{i} \telt$ a scalar pair $(\myqX_{\telt,\selt},\myqY_{\selt,\telt})$, at least one of which is nonzero and with edge product $\myqP_{\selt,\telt} = \myqX_{\telt,\selt}\myqY_{\selt,\telt}$. 
For any fixed $\telt \in L$ and $i \in \{1,\ldots,n-1\}$, the associated {\em crossing relation} is the equation 
\begin{equation}
\sum_{\selt: \selt \myarrow{i} \telt}\myqP_{\selt,\telt} - \sum_{\uelt: \telt \myarrow{i} \uelt}\myqP_{\telt,\uelt} = \rho_{i}(\telt)-\delta_{i}(\telt).
\end{equation}
Given a diamond \parbox{1.4cm}{\begin{center}
\setlength{\unitlength}{0.2cm}
\begin{picture}(5,3)
\put(2.5,0){\circle*{0.5}} \put(0.5,2){\circle*{0.5}}
\put(2.5,4){\circle*{0.5}} \put(4.5,2){\circle*{0.5}}
\put(0.5,2){\line(1,1){2}} \put(2.5,0){\line(-1,1){2}}
\put(4.5,2){\line(-1,1){2}} \put(2.5,0){\line(1,1){2}}
\put(1.25,0.55){\em \small j} \put(3.2,0.7){\em \small i}
\put(1.2,2.7){\em \small i} \put(3.25,2.55){\em \small j}
\put(3,-0.75){\footnotesize $\relt$} \put(5.25,1.75){\footnotesize $\telt$}
\put(3,4){\footnotesize $\uelt$} \put(-1,1.75){\footnotesize $\selt$}
\end{picture} \end{center}}, the associated {\em diamond relations} are 
\begin{equation}
\myqX_{\uelt,\selt}\myqY_{\telt,\uelt} = \myqY_{\relt,\selt}\myqX_{\telt,\relt} \hspace*{0.2in}\mbox{ and }\hspace*{0.2in} \myqX_{\uelt,\telt}\myqY_{\selt,\uelt} = \myqY_{\relt,\telt}\myqX_{\selt,\relt}.
\end{equation}
Call the set of all these diamond and crossing relations the {\em DC relations}. 
Next, let $V[L]$ be a complex vector space freely generated by $\{v_{\telt}\}_{\telt \in L}$. 
For all $\selt, \telt \in L$ and all $i \in \{1,\ldots,n-1\}$, define  
\begin{equation}
\myqx_{i}.v_{\selt} := \sum_{\uelt \in L}\myqX_{\uelt,\selt}v_{\uelt} \hspace*{0.2in}\mbox{ and }\hspace*{0.2in} \myqy_{i}.v_{\telt} := \sum_{\relt \in L}\myqY_{\relt,\telt}v_{\relt}.
\end{equation}
A version of the next result originally appeared in \cite{DonSupp}, but the statement below follows Theorems 10.8/10.9/10.10 of the manuscript \cite{DonDiamond}. 

\noindent 
{\bf \SetupProposition}\ \ {\sl With scalar pairs assigned to edges as above, $L$ is} $\myA_{n-1}${\sl -structured and all DC relations are satisfied if and only if (i) equations (3) above induce a well-defined action of} $\mathfrak{g}(\myA_{n-1})$ {\sl on $V[L]$ and (ii) $\{v_{\telt}\}_{\telt \in L}$ is a weight basis for the} $\mathfrak{g}(\myA_{n-1})${\sl -module $V[L]$ and (iii) the diamond-colored distributive lattice $L$ together with the assigned scalars is the representation diagram for this weight basis. In this case,} $\WGF(L;\myvarZ)$ {\sl is} $W(\myA_{n-1})$-{\sl invariant, and $V[L]$ is isomorphic to the direct sum $V(\lambda_{1}) \oplus \cdots \oplus V(\lambda_{l})$ of irreducible} $\mathfrak{g}(\myA_{n-1})${\sl -modules if and only if} $\WGF(L;\myvarZ) = \sum_{1 \leq k \leq l} \WGF(R(\lambda_{k});\myvarZ)$. 

We close this section with a brief discussion of some `extremal' properties enjoyed by some weight bases and their supporting graphs, cf.\ \cite{DonSupp}. 
A supporting graph $R$ for a given $\mathfrak{g}(\myA_{n-1})$-module $V$, and the weight basis $\mathscr{B} = \{v_{\telt}\}_{\telt \in R}$ associated with $R$, is {\em edge-minimal} if no subgraph of $R$ is the supporting graph for another weight basis of $V$. 
Let $R' := \{\telt'\, |\, \telt \in R\}$ be a set of symbols. 
We say $R$ (and $\mathscr{B}$) is {\em solitary} if, whenever $\mathscr{B}' = \{v_{\telt'}\}_{\telt' \in R'}$ is another weight basis for $V$ whose supporting graph $R'$ is isomorphic to $R$ as an edge-colored directed graph under the correspondence $\telt' \leftrightarrow \telt$, then there exist scalars $\{c_{\telt}\}_{\telt \in R'}$ such that $\velt_{\telt'} = c_{\telt}v_{\telt}$ for all $\telt  \in R$. 
A solitary weight basis can be thought of as being uniquely identified, in the above sense, by its supporting graph. 
We say two representation diagrams $R$ and $R'$ are {\em edge-product similar} if they are isomorphic as edge-colored directed graphs and if edge products are preserved under some such isomorphism. 
Notice that if $R$ is solitary, then any two representation diagrams whose underlying supporting graph is $R$ must be edge-product similar. 
The next result says how some of these algebraic and order-theoretic notions can interact. 

\noindent 
{\bf \SolitaryLatticeProposition}\ \ {\sl Suppose $L$ is a diamond-colored distributive lattice representation diagram for some} $\mathfrak{g}(\myA_{n-1})$-{\sl module $V$. 
Further, suppose the DC relations on $L$ uniquely determine its edge products, and suppose no edge product is equal to zero.  
Then $L$ is edge-minimal. 
Moreover, if $V$ is irreducible, then $L$ is solitary.} 

{\em Proof.} In \cite{DonDiamond}, see Theorem 10.10 parts (3) and (4) for proofs (respectively) of the solitarity and edge-minimality claims of the proposition statement.\hfill\QED

\newcommand{\ThreeThreeTab}[4]
{\setlength{\unitlength}{0.2cm}\begin{picture}(3.2,1)\put(0.1,-0.75){\line(0,1){1}} \put(1.1,-0.75){\line(0,1){2}} \put(2.1,-0.75){\line(0,1){2}} \put(3.1,0.25){\line(0,1){1}} \put(0.1,-0.75){\line(1,0){2}} \put(0.1,0.25){\line(1,0){3}} \put(1.1,1.25){\line(1,0){2}} \put(0.35,-0.5){\tiny #1} \put(1.35,0.5){\tiny #2} \put(1.35,-0.5){\tiny #3} \put(2.35,0.5){\tiny #4} \end{picture}}
\vspace*{0.1in}
\noindent {\bf \S \MainResultSection\ Symmetric Fibonaccian lattices are skew-tabular lattices.} 
The main result of this section, and of the paper, is a correspondence that, once observed, is rather elementary. 
Specifically, \FibGTTheorem\ records our observation that our symmetric Fibonaccian lattices are a special family of the skew-tabular lattices introduced in \cite{DD}.\footnote{The first- and second-named authors actually discovered the results of \cite{DD} in an effort to understand the representation-theoretic aspects of the symmetric Fibonaccian lattices considered in the present paper.  That is, this project is the progenitor \cite{DD}, and not vice-versa.} 

Take as given two integers $m$ and $n$ with $m \geq n \geq 2$ and two nonincreasing nonnegative-integer sequences $\mysmallP = (\mysmallP_{1},\mysmallP_{2},\ldots,\mysmallP_{m})$ and $\mysmallQ = (\mysmallQ_{1},\mysmallQ_{2},\ldots,\mysmallQ_{m})$ with $\mysmallP_{i} \geq \mysmallQ_{i}$ for all $i \in \{1,2,\ldots,m\}$. 
So, $\mysmallQ$ is an integer partition whose partition diagram (aka Ferrers diagram) is contained within the partition diagram for the integer partition $\mysmallP$. 
The cells of the partition $\mysmallP$ are indexed by (row,column) pairs, as with matrices. 
We keep this same indexing when we consider the {\em skew shape} $\mysmallP/\mysmallQ$ obtained by removing the cells of $\mysmallQ$ from $\mysmallP$. 
Say $\mysmallP$ and $\mysmallQ$ are {\em $n$-skew-compatible} if no column of the skew shape $\mysmallP/\mysmallQ$ has more than $n$ cells. 
A {\em semistandard $n$-tableau $T$ of shape} $\mysmallP/\mysmallQ$ is a collection of integers $(T_{i,j})_{(i,j) \in \mytinyP/\mytinyQ}$, such that, for all cells $(i,j) \in \mysmallP/\mysmallQ$, we have $T_{i,j} \in \{1,2,\ldots,n\}$, $T_{i,j} \leq T_{i,j+1}$ (if $(i,j+1) \in \mysmallP/\mysmallQ$), and $T_{i,j} < T_{i+1,j}$ (if $(i+1,j) \in \mysmallP/\mysmallQ$). 
The {\em skew-tabular lattice} $L_{\mytinyA_{n-1}}^{\mbox{\tiny skew}}(\mysmallP/\mysmallQ)$ is the diamond-colored distributive lattice wherein $S \rightarrow T$ is a covering relation if and only if there exists some cell $(a,b) \in \mysmallP/\mysmallQ$ such that $S_{p,q} = T_{p,q}$ for all cells $(p,q) \ne (a,b)$ while $S_{a,b} = T_{a,b}+1$; in this case we assign color $i := T_{a,b}$ to this directed edge and write $S \myarrow{i} T$. 
One can see that, for each $T \in L_{\mytinyA_{n-1}}^{\mbox{\tiny skew}}(\mysmallP/\mysmallQ)$ and $i \in \{1,\ldots,n-1\}$, we have $\mym_{i}(T) = \rho_{i}(T)-\delta_{i}(T) = \#_{i}(T)-\#_{i+1}(T)$, where $\#_{c}(X)$ denotes the number of cells of the semistandard $n$-tableau $X$ which contain the number $c$. 
For such $T$ we have, then, $\wt(T) = \sum_{i=1}^{n-1}\mym_{i}(T)\omega_{i} = \sum_{i=1}^{n-1}(\#_{i}(T)-\#_{i+1}(T))\omega_{i}$. 
Following \cite{DD}, we call $\vartheta_{\mytinyP/\mytinyQ} := \WGF(L_{\mytinyA_{n-1}}^{\mbox{\tiny skew}}(\mysmallP/\mysmallQ);\myvarZ)$ the {\em skew Schur function} for the shape $\mysmallP/\mysmallQ$.\footnote{Our $\vartheta_{_{\mytinyP/\mytinyQ}}$ differs mildly from the classical skew Schur function $\mys_{_{\mytinyP/\mytinyQ}}$, which is defined to be the following polynomial in $\mathbb{Z}[x_{1},x_{2},\ldots,x_{n}]$:  
\[\mys_{_{\mytinyP/\mytinyQ}} := \sum_{T \in L_{\mytinyA_{n-1}}^{\mbox{\tiny skew}}(\mytinyP/\mytinyQ)}x_{1}^{\#_1(T)}x_{2}^{\#_2(T)}\cdots x_{n}^{\#_n(T)}.\]
Observe that $\vartheta_{_{\mytinyP/\mytinyQ}}$ is the image of $\mys_{_{\mytinyP/\mytinyQ}}$ under the change of variables $x_{i} = z_{i-1}^{-1}z_{i}$ taking $z_{0} := 1 =: z_{n}$.}  
For an example, see \GTFig. 
When $\mysmallQ = (0,\ldots,0)$, the skew-tabular lattice is what we call a `classical GT lattice' $L_{\mytinyA_{n-1}}^{\mbox{\tiny GT}}(\mysmallP)$ whose weight-generating function is the {\em Schur function} $\vartheta_{\mytinyP}$.

A {\em right-to-left (RTL) sequence} in $\mysmallP/\mysmallQ$, or in a given tableau of shape $\mysmallP/\mysmallQ$, is a nonempty sequence of distinct cells that begins at the upper-rightmost cell and has the property that every other cell in the sequence is preceded by the cell to its right on the same row or, if there is no such cell, is preceded by the leftmost cell in the row above. 
A tableau $T$ from the skew-tabular lattice $L_{\mytinyA_{n-1}}^{\mbox{\tiny skew}}(\mysmallP/\mysmallQ)$ is {\em ballot-admissible}\footnote{The more common term is {\em Littlewood--Richardson (LR) tableau}; see Chapter 7 Appendix 1.3 of \cite{StanText2} for some discussion relating to such terminology.} if, for all $i \in \{1,\ldots,n-1\}$ and each RTL sequence in $T$, the number of cells with an $i+1$ entry does not exceed the number of cells with an entry of $i$.
Let $\mathcal{B}(\mysmallP/\mysmallQ)$ denote the set of ballot-admissible semistandard $n$-tableaux of shape $\mysmallP/\mysmallQ$. 
(For example, in \GTFig, we have $\mathcal{B}\left(\rule[-1.5mm]{-0.05mm}{4mm}\ThreeThreeTab{ }{ }{ }{ }\right) = \left\{\rule[-1.5mm]{-0.05mm}{4mm}\ThreeThreeTab{1}{1}{2}{1},\ThreeThreeTab{2}{1}{2}{1}\right\}$.)
For any $B \in \mathcal{B}(\mysmallP/\mysmallQ)$, we have $\mym_{i}(B) \geq 0$ for all $i \in \{1,\ldots,n-1\}$. 
In this case, we let $\mypart(B)$ be the integer partition $\left(\sum_{j=i}^{n-1} \mym_{i}(B)\right)_{i=1}^{n}$. 

Part (1) of the next result follows Theorem 6.4 of \cite{DD}.  
The weight generating function identity and $\mathfrak{g}(\myA_{n-1})$-module isomorphism claimed in Part (2) is a consequence of the well-known `Littlewood--Richardson Rule' (see, for example, Appendix 1 of Chapter 7 of \cite{StanText2}). 
A version of the Littlewood--Richardson Rule that is compatible with the notation and language of this paper is Theorem 5.2 of \cite{DD}. 

\noindent 
{\bf \DDTheorem\ (Donnelly--Dunkum)}\ \ {\sl Let} $L:=L_{\mytinyA_{n-1}}^{\mbox{\tiny skew}}(\mysmallP/\mysmallQ)$.  {\sl (1) Then $L$ is} $\myA_{n-1}$-{\sl structured.  Moreover, there exists a set of positive rational scalar pairs} $\left\{\rule[-1.2mm]{-0.05mm}{4mm}(\myqX_{T,S},\myqY_{S,T})\right\}_{S \rightarrow T \mbox{ in } L}$ {\sl that, when assigned to the edges of $L$, satisfy all DC relations.  (2) In particular, $L$ is a representation diagram for the} $\mathfrak{g}(\myA_{n-1})$-{\sl module $V[L]$ whose associated} $W(\myA_{n-1})$-{\sl invariant weight-generating function is}
\[\WGF(L;\myvarZ) = \vartheta_{\mytinyP/\mytinyQ} = \sum_{B \in \mathcal{B}(\mytinyP/\mytinyQ)}\vartheta_{\mysmallerpart(B)}\]
{\sl and which is therefore isomorphic to the following direct sum of irreducible} $\mathfrak{g}(\myA_{n-1})$-{\sl modules:} 
\[V[L] \cong \bigoplus_{B \in \mathcal{B}(\mytinyP/\mytinyQ)} V(\wt(B)).\]

Fix for the remainder of this section positive integers $n$ and $k$ with $n \geq 2$, and let $\mysmallP$ and $\mysmallQ$ be the partitions defined by the following weakly decreasing nonnegative integer $\left(\rule[-1.5mm]{0mm}{4.75mm}1+\lfloor\frac{k}{2}\rfloor (n-2)\right)$-tuples:
\begin{eqnarray*}
\mysmallP & := & \left(\rule[-1.5mm]{0mm}{4.75mm}k,\underbrace{k-1,\ldots,k-1}_{n-2},\underbrace{k-3,\ldots,k-3}_{n-2},\ldots,\underbrace{1+(k\, \mbox{\footnotesize mod}\, 2),\ldots,1+(k\, \mbox{\footnotesize mod}\, 2)}_{n-2}\right)\\
\mysmallQ & := & \left(\rule[-1.5mm]{0mm}{4.75mm}\underbrace{k-2,\ldots,k-2}_{n-2},\underbrace{k-4,\ldots,k-4}_{n-2},\ldots,\underbrace{k\, \mbox{\footnotesize mod}\, 2,\ldots,k\, \mbox{\footnotesize mod}\, 2}_{n-2},0\right).
\end{eqnarray*}
The skew shape $\mysmallR^{\mbox{\tiny Fib$(n,k)$}}_{\mbox{\tiny ribbon}} := \mysmallP/\mysmallQ$ is the {\em Fibonacci ribbon} consisting of those cells from the partition diagram of $\mysmallP$ that are bottommost in their column or rightmost in their row.  
This skew shape has exactly $1+\lfloor\frac{k}{2}\rfloor (n-2)$ nonempty rows.  
It also has exactly $k$ columns, and these column lengths toggle between $1$ and $n-1$, with $1$ as the length of the rightmost column. 
See \RibbonFigure. 

\begin{figure}[ht]
\begin{center}
\setlength{\unitlength}{0.35cm}
\begin{picture}(36,12)
\put(0,10.5){{\bf \RibbonFigure:} Fibonacci ribbons $\mysmallR^{\mbox{\tiny Fib$(n,k)$}}_{\mbox{\tiny ribbon}}$ for four different pairs $(n,k)$.}
\put(0,3.5){\begin{picture}(0,0)
\put(-1,3){\footnotesize $(3,3)$}
\put(0,0){\line(0,1){1}}
\put(0,0){\line(1,0){2}}
\put(1,0){\line(0,1){2}}
\put(0,1){\line(1,0){3}}
\put(2,0){\line(0,1){2}}
\put(1,2){\line(1,0){2}}
\put(3,1){\line(0,1){1}}
\end{picture}}
\put(10,0){\begin{picture}(0,0)
\put(-1,8.5){\footnotesize $(6,4)$}
\put(0,0){\line(0,1){5}}
\put(1,0){\line(0,1){5}}
\multiput(0,0)(0,1){6}{\line(1,0){1}}
\put(1,4){\line(1,0){1}}
\put(1,5){\line(1,0){1}}
\put(2,4){\line(0,1){5}}
\put(3,4){\line(0,1){5}}
\multiput(2,4)(0,1){6}{\line(1,0){1}}
\put(3,8){\line(1,0){1}}
\put(3,9){\line(1,0){1}}
\put(4,8){\line(0,1){1}}
\end{picture}}
\put(20,0.5){\begin{picture}(0,0)
\put(-1,6.5){\footnotesize $(5,4)$}
\put(0,0){\line(0,1){4}}
\put(1,0){\line(0,1){4}}
\multiput(0,0)(0,1){5}{\line(1,0){1}}
\put(1,3){\line(1,0){1}}
\put(1,4){\line(1,0){1}}
\put(2,3){\line(0,1){4}}
\put(3,3){\line(0,1){4}}
\multiput(2,3)(0,1){5}{\line(1,0){1}}
\put(3,6){\line(1,0){1}}
\put(3,7){\line(1,0){1}}
\put(4,6){\line(0,1){1}}
\end{picture}}
\put(30,1.5){\begin{picture}(0,0)
\put(-1,4.5){\footnotesize $(4,5)$}
\put(0,0){\line(0,1){1}}
\put(0,0){\line(1,0){1}}
\put(0,1){\line(1,0){1}}
\put(1,0){\line(0,1){3}}
\put(2,0){\line(0,1){3}}
\multiput(1,0)(0,1){4}{\line(1,0){1}}
\put(2,2){\line(1,0){1}}
\put(2,3){\line(1,0){1}}
\put(3,2){\line(0,1){3}}
\put(4,2){\line(0,1){3}}
\multiput(3,2)(0,1){4}{\line(1,0){1}}
\put(4,4){\line(1,0){1}}
\put(4,5){\line(1,0){1}}
\put(5,4){\line(0,1){1}}
\end{picture}}
\end{picture}
\end{center}
\end{figure}

We now present an edge-and-edge-color-preserving poset isomorphism from the symmetric Fibonaccian lattice $L^{\mbox{\tiny Fib}}(n,k)$ to the skew-tabular lattice $L_{\mytinyA_{n-1}}^{\mbox{\tiny skew}}\!\!\left(\rule[-1.5mm]{0mm}{4.75mm}\mysmallR^{\mbox{\tiny Fib$(n,k)$}}_{\mbox{\tiny ribbon}}\right)$. 

\noindent
{\bf \FibGTTheorem}\ \ {\sl For $X = (X_{1},\ldots,X_{k})$ in} $L^{\mbox{\tiny Fib}}(n,k)$, {\sl let $Y$ be the tableau of shape} $\mysmallR^{\mbox{\tiny Fib$(n,k)$}}_{\mbox{\tiny ribbon}}$ {\sl whose $j^{\mbox{\tiny th}}$ column, starting from the left, has the set of entries $\{X_{j}-(j-1)n\}$ if $k+1-j$ is odd and $\{1,2,\ldots,n\} \setminus \{n+1-(X_{j}-(j-1)n)\}$ if $k+1-j$ is even, with these entries strictly increasing from top to bottom within the column. 
Then $Y$ is in} $L_{\mytinyA_{n-1}}^{\mbox{\tiny skew}}\!\!\left(\rule[-1.5mm]{0mm}{4.75mm}\mysmallR^{\mbox{\tiny Fib$(n,k)$}}_{\mbox{\tiny ribbon}}\right)$. 
{\sl Now define} $\varphi: L^{\mbox{\tiny Fib}}(n,k) \longrightarrow L_{\mytinyA_{n-1}}^{\mbox{\tiny skew}}\!\!\left(\rule[-1.5mm]{0mm}{4.75mm}\mysmallR^{\mbox{\tiny Fib$(n,k)$}}_{\mbox{\tiny ribbon}}\right)$ {\sl by the rule $\varphi(X) = Y$. 
Then $\varphi$ is an edge-and-edge-color-preserving poset isomorphism.} 

For an illustration of \FibGTTheorem, compare \IntroFibFig\ and \GTFig.

\begin{figure}[t]
\begin{center}
{\bf \GTFig:}  
The lattice $L := L_{\mytinyA_{2}}^{\mbox{\tiny skew}}\!\!\left(\rule[-1.5mm]{0mm}{4.75mm}\mysmallR^{\mbox{\tiny Fib$(3,3)$}}_{\mbox{\tiny ribbon}}\right) \cong L^{\mbox{\tiny Fib}}(3,3)$. For $L$ we have: 
\begin{eqnarray*} 
\RGF(L;q) & = & 1+3q+4q^{2}+5q^{3}+4q^{4}+3q^{5}+q^{6}\\ 
 & = & \left[q^{3\cdot{2}-\rho\left(\rule[-1.5mm]{-0.05mm}{4mm}\ThreeThreeTab{1}{1}{2}{1}\right)}\right]\cdot\frac{(1-q^{3})(1-q^{5})(1-q^{2})}{(1-q)(1-q^{2})(1-q)} + \left[q^{3\cdot{2}-\rho\left(\rule[-1.5mm]{-0.05mm}{4mm}\ThreeThreeTab{2}{1}{2}{1}\right)}\right]\cdot\frac{(1-q)(1-q^{4})(1-q^{3})}{(1-q)(1-q^{2})(1-q)}\\
\WGF(L;\myvarZ) & = & z_{1}^{2}z_{2} + 2z_{2}^{2} + z_{1}^{3}z_{2}^{-1} + 3z_{1} + z_{1}^{-2}z_{2}^{3}+2z_{1}^{2}z_{2}^{-2}\\ 
 & & \hspace*{2in}+ 3z_{1}^{-1}z_{2} + z_{1}^{-3}z_{2}^{2} + 3z_{2}^{-1} + z_{1}z_{2}^{-3} + 2z_{1}^{-2} + z_{1}^{-1}z_{2}^{-2}
\end{eqnarray*}
\setlength{\unitlength}{1.5cm}
\begin{picture}(4,6.5)
\put(0,0){\qbezier(2,6)(3.25,5)(4.5,4)}
\put(0,0){\qbezier(2,6)(1.5,5.5)(1,5)}
\put(0,0){\qbezier(2,6)(2,5.5)(2,5)}
\put(0,0){\qbezier(3.25,5)(2.75,4.5)(2.25,4)}
\put(0,0){\qbezier(3.25,5)(3.25,4)(3.25,3)}
\put(0,0){\qbezier(4.5,4)(4,3.5)(3.5,3)}
\put(0,0){\qbezier(4.5,4)(4.5,3)(4.5,2)}
\put(0,0){\qbezier(1,5)(2.25,4)(3.5,3)}
\put(0,0){\qbezier(1,5)(1,4.5)(1,4)}
\put(0,0){\qbezier(2,5)(3.25,4)(4.5,3)}
\put(0,0){\qbezier(2,5)(1,4)(0,3)}
\put(0,0){\qbezier(1,4)(2.25,3)(3.5,2)}
\put(0,0){\qbezier(2.25,4)(2.25,3)(2.25,2)}
\put(0,0){\qbezier(3.5,3)(3.5,2)(3.5,1)}
\put(0,0){\qbezier(3.25,4)(2.25,3)(1.25,2)}
\put(0,0){\qbezier(4.5,3)(3.5,2)(2.5,1)}
\put(0,0){\qbezier(0,3)(1.25,2)(2.5,1)}
\put(0,0){\qbezier(1.25,2)(1.25,1.5)(1.25,1)}
\put(0,0){\qbezier(2.5,1)(2.5,0.5)(2.5,0)}
\put(0,0){\qbezier(3.25,3)(2.25,2)(1.25,1)}
\put(0,0){\qbezier(3.25,3)(3.875,2.5)(4.5,2)}
\put(0,0){\qbezier(4.5,2)(3.5,1)(2.5,0)}
\put(0,0){\qbezier(2.25,2)(2.875,1.5)(3.5,1)}
\put(0,0){\qbezier(1.25,1)(1.875,0.5)(2.5,0)}
\put(2,6){\VertexTableauFibTwo{1}{1}{2}{1}{-0.6}{-0.05}}
\put(3.25,5){\VertexTableauFibTwo{1}{1}{2}{2}{0.125}{0}}
\put(4.5,4){\VertexTableauFibTwo{1}{1}{2}{3}{0.15}{-0.05}}
\put(2,5){\VertexTableauFibTwo{1}{1}{3}{1}{-0.575}{-0.075}}
\put(3.25,4){\VertexTableauFibTwo{1}{1}{3}{2}{0.1}{-0.025}}
\put(4.5,3){\VertexTableauFibTwo{1}{1}{3}{3}{0.15}{-0.05}}
\put(3.25,3){\VertexTableauFibTwo{1}{2}{3}{2}{-0.625}{-0.2}}
\put(4.5,2){\VertexTableauFibTwo{1}{2}{3}{3}{0.15}{-0.05}}
\put(1,5){\VertexTableauFibTwo{2}{1}{2}{1}{-0.6}{-0.05}}
\put(2.25,4){\VertexTableauFibTwo{2}{1}{2}{2}{-0.58}{-0.4}}
\put(3.5,3){\VertexTableauFibTwo{2}{1}{2}{3}{0.12}{-0.1}}
\put(1,4){\VertexTableauFibTwo{2}{1}{3}{1}{-0.6}{-0.05}}
\put(2.25,3){\VertexTableauFibTwo{2}{1}{3}{2}{-0.725}{-0.25}}
\put(3.5,2){\VertexTableauFibTwo{2}{1}{3}{3}{0.125}{-0.2}}
\put(2.25,2){\VertexTableauFibTwo{2}{2}{3}{2}{-0.625}{-0.05}}
\put(3.5,1){\VertexTableauFibTwo{2}{2}{3}{3}{0.15}{-0.1}}
\put(0,3){\VertexTableauFibTwo{3}{1}{3}{1}{-0.6}{-0.05}}
\put(1.25,2){\VertexTableauFibTwo{3}{1}{3}{2}{-0.625}{-0.3}}
\put(2.5,1){\VertexTableauFibTwo{3}{1}{3}{3}{-0.625}{-0.3}}
\put(1.25,1){\VertexTableauFibTwo{3}{2}{3}{2}{-0.625}{-0.05}}
\put(2.5,0){\VertexTableauFibTwo{3}{2}{3}{3}{0.15}{-0.1}}
\put(1,5){\NEEdgeLabelForLatticeI{{\em 1}}}
\put(2,3.75){\NEEdgeLabelForLatticeI{{\em 1}}}
\put(2.1,2.85){\NEEdgeLabelForLatticeI{{\em 1}}}
\put(2.1,1.85){\NEEdgeLabelForLatticeI{{\em 1}}}
\put(3.75,3.25){\NEEdgeLabelForLatticeI{{\em 1}}}
\put(3.75,2.25){\NEEdgeLabelForLatticeI{{\em 1}}}
\put(3.6,1.1){\NEEdgeLabelForLatticeI{{\em 1}}}
\put(1.3,4.3){\NEEdgeLabelForLatticeI{{\em 1}}}
\put(2,5){\VerticalEdgeLabelForLatticeI{{\em 2}}}
\put(1,4){\VerticalEdgeLabelForLatticeI{{\em 2}}}
\put(3.25,4){\VerticalEdgeLabelForLatticeI{{\em 2}}}
\put(2.25,3){\VerticalEdgeLabelForLatticeI{{\em 2}}}
\put(4.5,3){\VerticalEdgeLabelForLatticeI{{\em 2}}}
\put(3.5,2){\VerticalEdgeLabelForLatticeI{{\em 2}}}
\put(3.25,3){\VerticalEdgeLabelForLatticeI{{\em 1}}}
\put(2.25,2){\VerticalEdgeLabelForLatticeI{{\em 1}}}
\put(1.25,1){\VerticalEdgeLabelForLatticeI{{\em 1}}}
\put(4.5,2){\VerticalEdgeLabelForLatticeI{{\em 1}}}
\put(3.5,1){\VerticalEdgeLabelForLatticeI{{\em 1}}}
\put(2.5,0){\VerticalEdgeLabelForLatticeI{{\em 1}}}
\put(3.15,5){\NWEdgeLabelForLatticeI{{\em 1}}}
\put(2.4,3.8){\NWEdgeLabelForLatticeI{{\em 1}}}
\put(2.9,4.2){\NWEdgeLabelForLatticeI{{\em 1}}}
\put(2.025,3.1){\NWEdgeLabelForLatticeI{{\em 1}}}
\put(1.025,2.1){\NWEdgeLabelForLatticeI{{\em 1}}}
\put(4.525,3.9){\NWEdgeLabelForLatticeI{{\em 2}}}
\put(4.15,3.2){\NWEdgeLabelForLatticeI{{\em 2}}}
\put(4.65,1.8){\NWEdgeLabelForLatticeI{{\em 2}}}
\put(3.525,1.9){\NWEdgeLabelForLatticeI{{\em 2}}}
\put(2.1,3.3){\NEEdgeLabelForLatticeI{{\em 2}}}
\put(2.1,1.3){\NEEdgeLabelForLatticeI{{\em 2}}}
\put(1.6,0.9){\NEEdgeLabelForLatticeI{{\em 2}}}
\put(1.225,0.2){\NEEdgeLabelForLatticeI{{\em 2}}}
\put(0.1,3.1){\NEEdgeLabelForLatticeI{{\em 2}}}
\put(1.25,2){\NEEdgeLabelForLatticeI{{\em 2}}}
\put(2.45,1.3){\NWEdgeLabelForLatticeI{{\em 2}}}
\put(2.3,0.8){\NEEdgeLabelForLatticeI{{\em 2}}}
\put(3.5,0.1){\NWEdgeLabelForLatticeI{{\em 2}}}
\end{picture}
\end{center}
\end{figure}

{\em Proof of \FibGTTheorem.} First, we verify that $Y$ is in $L_{\mytinyA_{n-1}}^{\mbox{\tiny skew}}\!\!\left(\rule[-1.5mm]{0mm}{4.75mm}\mysmallR^{\mbox{\tiny Fib$(n,k)$}}_{\mbox{\tiny ribbon}}\right)$. 
Suppose columns $Y_{j}$ and $Y_{j+1}$ are obtained from $X_{j}$ and $X_{j+1}$ respectively. 
We consider two cases: $k+1-j$ is odd and $k+1-j$ is even. 
When $k+1-j$ is odd, then column $Y_{j}$ is the singleton $X_{j}-(j-1)n$, and column $Y_{j+1}$ has entries $1$, $2$, $\ldots$, $n-(X_{j+1}-jn)$, $n+2-(X_{j+1}-jn)$, $\ldots$, $n$ in order reading from the top of the column. 
We must check that the bottommost entry $b$ of $Y_{j+1}$ is no smaller than $a := X_{j}-(j-1)n$. 
Since $(j-1)n+1 \leq X_{j} \leq jn$, then $1 \leq a \leq n$. 
Also, $b$ must either be $n-1$ or $n$. 
So, it suffices to show that if $a=n$, then $b = n$. 
Well, if $a=n$, then $X_{j}-(j-1)n=n$, i.e.\ $X_{j} = jn$. 
Therefore $jn+1 < X_{j+1}$, from which it follows that $n+1-(X_{j+1}-jn) < n$. 
This latter fact forces us to have $b=n$.
Similarly, when $k+1-j$ is even, we must check that the topmost entry of column $Y_{j}$ is no larger than the singleton entry $X_{j+1}-jn$ of column $Y_{j+1}$. 
This is easily verified with reasoning that is similar to the preceding case. 
Since the procedure defining $\varphi$ reverses, then $\varphi$ is a bijection. 
It is routine to check that $S \myarrow{i} T$ if and only if $\varphi(S) \myarrow{i} \varphi(T)$, which concludes the argument.
\hfill\QED

\noindent
{\bf \FibGTCorollary}\ \ {\sl The symmetric Fibonaccian lattice} $L^{\mbox{\tiny Fib}}(n,k)$ {\sl is the supporting graph for a weight basis for a representation of} $\mathfrak{g}(\myA_{n-1}) \cong \mathfrak{sl}(n,\mathbb{C})$ {\sl whose associated} $W(\myA_{n-1})$-{\sl invariant weight-generating function is the skew Schur function} {\Large $\vartheta_{\mysmallR^{\mbox{\tiny Fib$(n,k)$}}_{\mbox{\tiny ribbon}}}$}. 

{\em Proof.} Follows from \FibGTTheorem\ together with \DDTheorem.\hfill\QED

In view of \FibGTTheorem, it now follows from Corollary 6.6 of \cite{DD} that $L^{\mbox{\tiny Fib}}(n,k)$ is rank symmetric, rank unimodal, and strongly Sperner. 
In fact, a much stronger statement holds. 
In \S \IntroSection, we noted that, disregarding edge/vertex colors, the rank symmetric distributive lattice $L^{\mbox{\tiny Fib}}(n,k)$ is isomorphic to the distributive lattice $\mathbf{J}(Q(n-1,k(n-1)))$ of order ideals from Gansner's `up-down' poset $Q(n-1,k(n-1))$. 
A consequence of Theorem 2 from \cite{Gan} is that $\mathbf{J}(Q(n-1,k(n-1)))$ is a nested chain order. 
This nested chain order property together with rank symmetry implies that:  

\noindent
{\bf \SCDTheorem\ (Gansner)}\ \ {\sl The symmetric Fibonaccian lattice $L^{\mbox{\tiny Fib}}(n,k)$ has a symmetric chain decomposition.}

\vspace*{0.1in}
\noindent {\bf \S \RGFSection\ Rank generating functions.} 
Having realized symmetric Fibonaccian lattices as skew-tabular lattices, we can use results from \cite{DD} to better understand their rank generating functions. 
To do so, we identify a distinguished subset of Fibonaccian strings from $L^{\mbox{\tiny Fib}}(n,k)$ via the isomorphism $\varphi: L^{\mbox{\tiny Fib}}(n,k) \longrightarrow L_{\mytinyA_{n-1}}^{\mbox{\tiny skew}}\!\!\left(\rule[-1.5mm]{0mm}{4.75mm}\mysmallR^{\mbox{\tiny Fib$(n,k)$}}_{\mbox{\tiny ribbon}}\right)$ from \FibGTTheorem. 
Let $\mathcal{B}^{\mbox{\tiny Fib}}(n,k) := \{\varphi^{-1}(B) \in L^{\mbox{\tiny Fib}}(n,k)\, |\, B \in \mathcal{B}\!\left(\rule[-1.5mm]{0mm}{4.75mm}\mysmallR^{\mbox{\tiny Fib$(n,k)$}}_{\mbox{\tiny ribbon}}\right)\}$, so $\mathcal{B}^{\mbox{\tiny Fib}}(n,k)$ is the set of ballot-admissible Fibonaccian strings in $L^{\mbox{\tiny Fib}}(n,k)$. 

Before we state the main result of this section, we introduce some notation. For any positive integer $\ell$, let $[\ell]_{q}$ denote the `$q$-integer' $(1-q^{\ell})/(1-q)$.

\noindent
{\bf \RGFTheorem}\ \ {\sl Let $n$ and $k$ be positive integers, with $n \geq 2$. 
Set $\lambda_{i}^{j}(X) := \sum_{l=i}^{j}(\rho_{l}(X)-\delta_{l}(X))$ for any $X$ in} $L^{\mbox{\tiny Fib}}(n,k)$.  {\sl Then}
\[\RGF(L^{\mbox{\tiny Fib}}(n,k);q) = \sum_{\mbox{\scriptsize $B \in \mathcal{B}^{\mbox{\tiny Fib}}(n,k)$}} \left(q^{k(n-1)-\rho(B)}\prod_{i=1}^{n-1}\prod_{j=i}^{n-1}\frac{[\lambda_{i}^{j}(B)+j+1-i]_{q}}{[j+1-i]_{q}}\right) = {A}^{(n)}_{k}(q)\]
{\sl is a symmetric and unimodal polynomial of degree $k(n-1)$ in the variable $q$.} 

{\em Proof.} The first claimed equality concerning $\RGF(L^{\mbox{\tiny Fib}}(n,k);q)$ follows from Theorem 5.2 of \cite{DD}, via the isomorphism $\varphi$ of \FibGTTheorem. 
The equality $\RGF(L^{\mbox{\tiny Fib}}(n,k);q) = {A}^{(n)}_{k}(q)$ follows from \RankProposition. 
Our polynomial has degree $k(n-1)$ since $L^{\mbox{\tiny Fib}}(n,k)$ has length $k(n-1)$. 
Of course, symmetry of $\RGF(L^{\mbox{\tiny Fib}}(n,k);q)$ follows from the fact that, when we disregard edge colors, this lattice is self-dual. 
Unimodality of the rank generating function is a consequence of Corollary 6.6 of \cite{DD}, since $L^{\mbox{\tiny Fib}}(n,k)$ has been shown in \FibGTTheorem\ to be a skew-tabular lattice. 
Alternatively, rank unimodality is a consequence of Gansner's result (\SCDTheorem\ above) that $L^{\mbox{\tiny Fib}}(n,k)$ has a symmetric chain decomposition.\hfill\QED

\vspace*{0.1in}
\noindent {\bf \S \EndSection\ Connections, conjectures, and open problems.} 
Following \cite{Malone}, we connect the preceding work with some enumerative and algebraic problems. 
We begin with the following observation. 

\noindent 
{\bf \FibSolitaryProposition}\ \ {\sl Let} $L := L^{\mbox{\tiny Fib}}(n,k)$ {\sl be a symmetric Fibonaccian lattice. 
For $k \in \{1,2,3\}$ and any integer $n$ with $n \geq 3$, the DC relations on} $L$ {\sl uniquely determine its edge products.  
In particular, as a supporting graph for the associated} $\mathfrak{g}(\myA_{n-1})$-{\sl module, $L$ is edge-minimal and, when $k \in \{1,2\}$, solitary. 
For all positive integers $k$,} $L^{\mbox{\tiny Fib}}(2,k)$ {\sl is solitary and edge-minimal.} 

A sketch of the proof is given below. 
For reference, the symmetric Fibonaccian lattice of \LatticeFigList\ is presented again in \FibSolitaryFig\ with the unique edge product identified for each edge. 
We have verified by hand that each of $L^{\mbox{\tiny Fib}}(3,k)$ is solitary when $k\in \{3,4\}$.  
An instructive exercise is to verify that the unique edge products presented in \FibSolitaryFig\ are correct.
 
The results of \FibSolitaryProposition\ suggest the following question: 

\noindent
{\bf \FibSolitaryOpen}\ \ For which $n$ and $k$ is $L^{\mbox{\tiny Fib}}(n,k)$ edge-minimal and/or solitary? 

\begin{figure}[t]
\begin{center}
{\bf \FibSolitaryFig:}  {\small In this depiction of $L_{\mytinyA_{2}}^{\mbox{\tiny skew}}\!\!\left(\rule[-1.5mm]{0mm}{4.75mm}\mysmallR^{\mbox{\tiny Fib$(3,3)$}}_{\mbox{\tiny ribbon}}\right) \cong L^{\mbox{\tiny Fib}}(3,3)$, each circled number is the edge product $\myqP_{\selt,\telt}$ associated with an edge $\selt \rightarrow \telt$; these edge products are uniquely determined by the DC relations.}

\setlength{\unitlength}{1.5cm}
\begin{picture}(4,6.5)
\put(0,0){\qbezier(2,6)(3.25,5)(4.5,4)}
\put(0,0){\qbezier(2,6)(1.5,5.5)(1,5)}
\put(0,0){\qbezier(2,6)(2,5.5)(2,5)}
\put(0,0){\qbezier(3.25,5)(2.75,4.5)(2.25,4)}
\put(0,0){\qbezier(3.25,5)(3.25,4)(3.25,3)}
\put(0,0){\qbezier(4.5,4)(4,3.5)(3.5,3)}
\put(0,0){\qbezier(4.5,4)(4.5,3)(4.5,2)}
\put(0,0){\qbezier(1,5)(2.25,4)(3.5,3)}
\put(0,0){\qbezier(1,5)(1,4.5)(1,4)}
\put(0,0){\qbezier(2,5)(3.25,4)(4.5,3)}
\put(0,0){\qbezier(2,5)(1,4)(0,3)}
\put(0,0){\qbezier(1,4)(2.25,3)(3.5,2)}
\put(0,0){\qbezier(2.25,4)(2.25,3)(2.25,2)}
\put(0,0){\qbezier(3.5,3)(3.5,2)(3.5,1)}
\put(0,0){\qbezier(3.25,4)(2.25,3)(1.25,2)}
\put(0,0){\qbezier(4.5,3)(3.5,2)(2.5,1)}
\put(0,0){\qbezier(0,3)(1.25,2)(2.5,1)}
\put(0,0){\qbezier(1.25,2)(1.25,1.5)(1.25,1)}
\put(0,0){\qbezier(2.5,1)(2.5,0.5)(2.5,0)}
\put(0,0){\qbezier(3.25,3)(2.25,2)(1.25,1)}
\put(0,0){\qbezier(3.25,3)(3.875,2.5)(4.5,2)}
\put(0,0){\qbezier(4.5,2)(3.5,1)(2.5,0)}
\put(0,0){\qbezier(2.25,2)(2.875,1.5)(3.5,1)}
\put(0,0){\qbezier(1.25,1)(1.875,0.5)(2.5,0)}
\put(2,6){\VertexBlankFib}
\put(3.25,5){\VertexBlankFib}
\put(4.5,4){\VertexBlankFib}
\put(2,5){\VertexBlankFib}
\put(3.25,4){\VertexBlankFib}
\put(4.5,3){\VertexBlankFib}
\put(3.25,3){\VertexBlankFib}
\put(4.5,2){\VertexBlankFib}
\put(1,5){\VertexBlankFib}
\put(2.25,4){\VertexBlankFib}
\put(3.5,3){\VertexBlankFib}
\put(1,4){\VertexBlankFib}
\put(2.25,3){\VertexBlankFib}
\put(3.5,2){\VertexBlankFib}
\put(2.25,2){\VertexBlankFib}
\put(3.5,1){\VertexBlankFib}
\put(0,3){\VertexBlankFib}
\put(1.25,2){\VertexBlankFib}
\put(2.5,1){\VertexBlankFib}
\put(1.25,1){\VertexBlankFib}
\put(2.5,0){\VertexBlankFib}
\put(1.1,5.9){\line(1,-1){0.35}}
\put(0.8,6){\tiny $2/3$}
\put(0.92,6.05){\circle{0.4}}
\put(0,0){\qbezier(1.1,5.9)(1.45,5.15)(1.8,4.4)}
\put(2.87,6){\scriptsize $1$}
\put(2.92,6.05){\circle{0.4}}
\put(0,0){\qbezier(2.72,5.94)(2.41,5.74)(2.1,5.54)}
\put(0.15,3.95){\line(1,-1){0.35}}
\put(-0.08,4.05){\scriptsize $2$}
\put(-0.03,4.10){\circle{0.4}}
\put(0,0){\qbezier(0.18,4.22)(0.55,4.36)(0.92,4.5)}
\put(0.18,2.08){\line(1,1){0.35}}
\put(-0.03,1.85){\scriptsize $2$}
\put(0.02,1.90){\circle{0.4}}
\put(0,0){\qbezier(0.17,1.74)(0.665,1.62)(1.16,1.5)}
\put(1.2,0.3){\line(1,1){0.375}}
\put(0.9,0.1){\tiny $2/3$}
\put(1.02,0.15){\circle{0.4}}
\put(0,0){\qbezier(1.2,0.3)(1.575,0.975)(1.95,1.65)}
\put(1.73,5.13){\line(0,-1){0.35}}
\put(1.6,5.3){\tiny $1/3$}
\put(1.72,5.35){\circle{0.4}}
\put(2.33,5.13){\line(0,-1){0.35}}
\put(2.2,5.3){\tiny $8/3$}
\put(2.32,5.35){\circle{0.4}}
\put(3.4,5.8){\tiny $4/3$}
\put(3.52,5.85){\circle{0.4}}
\put(0,0){\qbezier(3.41,5.66)(2.965,4.99)(2.52,4.32)}
\put(0,0){\qbezier(3.41,5.66)(3.095,5.555)(2.78,5.45)}
\put(3.8,5.4){\tiny $1/2$}
\put(3.92,5.45){\circle{0.4}}
\put(0,0){\qbezier(3.81,5.26)(3.765,4.52)(3.72,3.78)}
\put(0,0){\qbezier(3.81,5.26)(3.59,4.905)(3.37,4.55)}
\put(3.55,3.33){\scriptsize $3$}
\put(3.6,3.38){\circle{0.3}}
\put(0,0){\qbezier(3.425,3.34)(3.3725,3.38)(3.32,3.42)}
\put(1.77,3.8){\scriptsize $2$}
\put(1.82,3.85){\circle{0.4}}
\put(0,0){\qbezier(2.0,3.75)(2.24,3.75)(2.48,3.75)}
\put(0.9,2.95){\tiny $4/3$}
\put(1.02,3.00){\circle{0.4}}
\put(0,0){\qbezier(1.08,3.21)(1.29,3.355)(1.5,3.5)}
\put(0,0){\qbezier(1.13,2.82)(1.4,2.685)(1.67,2.55)}
\put(1.47,2.95){\scriptsize $1$}
\put(1.52,3.00){\circle{0.4}}
\put(0,0){\qbezier(1.71,3.13)(1.955,3.29)(2.2,3.45)}
\put(0,0){\qbezier(1.68,2.82)(1.94,2.685)(2.2,2.55)}
\put(2.8,2.95){\tiny $2/3$}
\put(2.92,3.00){\circle{0.4}}
\put(0,0){\qbezier(2.72,3)(2.66,3.11)(2.6,3.22)}
\put(0,0){\qbezier(2.72,3)(2.67,2.925)(2.62,2.85)}
\put(3.87,2.95){\scriptsize $1$}
\put(3.92,3.00){\circle{0.4}}
\put(0,0){\qbezier(4,3.2)(4.1,3.4)(4.2,3.6)}
\put(0,0){\qbezier(4,2.8)(4.1,2.585)(4.2,2.37)}
\put(4.79,2.95){\tiny $3/2$}
\put(4.92,3.00){\circle{0.4}}
\put(0,0){\qbezier(4.83,3.2)(4.715,3.3)(4.6,3.4)}
\put(0,0){\qbezier(4.83,3.2)(4.475,3.725)(4.12,4.25)}
\put(0,0){\qbezier(4.8,2.8)(4.7,2.685)(4.6,2.57)}
\put(0,0){\qbezier(4.8,2.8)(4.46,2.235)(4.12,1.67)}
\put(3.1,2.55){\scriptsize $3$}
\put(3.15,2.60){\circle{0.3}}
\put(0,0){\qbezier(3.28,2.48)(3.365,2.46)(3.45,2.44)}
\put(3.3,0){\scriptsize $1$}
\put(3.35,0.05){\circle{0.4}}
\put(0,0){\qbezier(3.13,0.13)(2.865,0.285)(2.6,0.44)}
\put(3.75,0.45){\tiny $4/3$}
\put(3.87,0.5){\circle{0.4}}
\put(0,0){\qbezier(3.66,0.55)(3.14,1.1)(2.62,1.65)}
\put(0,0){\qbezier(3.66,0.55)(3.385,0.55)(3.11,0.55)}
\put(4.05,0.95){\tiny $1/2$}
\put(4.17,1.0){\circle{0.4}}
\put(0,0){\qbezier(4.03,1.16)(3.865,1.63)(3.7,2.1)}
\put(0,0){\qbezier(4.03,1.16)(3.805,1.28)(3.58,1.4)}
\put(2.08,0.9){\line(0,1){0.38}}
\put(1.95,0.65){\tiny $1/3$}
\put(2.07,0.7){\circle{0.4}}
\put(2.78,0.9){\line(0,1){0.29}}
\put(2.65,0.65){\tiny $8/3$}
\put(2.77,0.7){\circle{0.4}}
\put(1.77,2.1){\scriptsize $2$}
\put(1.82,2.15){\circle{0.4}}
\put(0,0){\qbezier(2.03,2.18)(2.26,2.235)(2.49,2.29)}
\put(1,5){\NEEdgeLabelForLatticeI{{\em 1}}}
\put(2,3.75){\NEEdgeLabelForLatticeI{{\em 1}}}
\put(2.1,2.85){\NEEdgeLabelForLatticeI{{\em 1}}}
\put(2.1,1.85){\NEEdgeLabelForLatticeI{{\em 1}}}
\put(3.75,3.25){\NEEdgeLabelForLatticeI{{\em 1}}}
\put(3.25,1.75){\NEEdgeLabelForLatticeI{{\em 1}}}
\put(3.6,1.1){\NEEdgeLabelForLatticeI{{\em 1}}}
\put(1.3,4.3){\NEEdgeLabelForLatticeI{{\em 1}}}
\put(2,5){\VerticalEdgeLabelForLatticeI{{\em 2}}}
\put(1,4){\VerticalEdgeLabelForLatticeI{{\em 2}}}
\put(3.25,4){\VerticalEdgeLabelForLatticeI{{\em 2}}}
\put(2.25,3){\VerticalEdgeLabelForLatticeI{{\em 2}}}
\put(4.5,3){\VerticalEdgeLabelForLatticeI{{\em 2}}}
\put(3.5,2){\VerticalEdgeLabelForLatticeI{{\em 2}}}
\put(3.25,3){\VerticalEdgeLabelForLatticeI{{\em 1}}}
\put(2.25,2){\VerticalEdgeLabelForLatticeI{{\em 1}}}
\put(1.25,1){\VerticalEdgeLabelForLatticeI{{\em 1}}}
\put(4.5,2){\VerticalEdgeLabelForLatticeI{{\em 1}}}
\put(3.5,1){\VerticalEdgeLabelForLatticeI{{\em 1}}}
\put(2.5,0){\VerticalEdgeLabelForLatticeI{{\em 1}}}
\put(3.15,5){\NWEdgeLabelForLatticeI{{\em 1}}}
\put(2.4,3.8){\NWEdgeLabelForLatticeI{{\em 1}}}
\put(2.9,4.2){\NWEdgeLabelForLatticeI{{\em 1}}}
\put(2.025,3.1){\NWEdgeLabelForLatticeI{{\em 1}}}
\put(1.025,2.1){\NWEdgeLabelForLatticeI{{\em 1}}}
\put(4.525,3.9){\NWEdgeLabelForLatticeI{{\em 2}}}
\put(4.15,3.2){\NWEdgeLabelForLatticeI{{\em 2}}}
\put(4.65,1.8){\NWEdgeLabelForLatticeI{{\em 2}}}
\put(3.025,2.3){\NWEdgeLabelForLatticeI{{\em 2}}}
\put(2.1,3.3){\NEEdgeLabelForLatticeI{{\em 2}}}
\put(2.1,1.3){\NEEdgeLabelForLatticeI{{\em 2}}}
\put(1.6,0.9){\NEEdgeLabelForLatticeI{{\em 2}}}
\put(1.225,0.2){\NEEdgeLabelForLatticeI{{\em 2}}}
\put(0.1,3.1){\NEEdgeLabelForLatticeI{{\em 2}}}
\put(1.25,2){\NEEdgeLabelForLatticeI{{\em 2}}}
\put(2.45,1.3){\NWEdgeLabelForLatticeI{{\em 2}}}
\put(2.3,0.8){\NEEdgeLabelForLatticeI{{\em 2}}}
\put(3.5,0.1){\NWEdgeLabelForLatticeI{{\em 2}}}
\end{picture}
\end{center}

\

\vspace*{-0.5in}
\end{figure}

{\em Sketch of proof of \FibSolitaryProposition.} 
When $(n,k) \in \{(2,m)\}_{m \geq 1} \cup \{(m,1)\}_{m \geq 3} \cup \{(m,2)\}_{m \geq 3}$, the Fibonacci ribbon $\mysmallR^{\mbox{\tiny Fib$(n,k)$}}_{\mbox{\tiny ribbon}}$ is non-skew. 
By Theorem 8.3.2 of \cite{DD}, it follows that  $L^{\mbox{\tiny Fib}}(n,k) \cong L_{\mytinyA_{n-1}}^{\mbox{\tiny skew}}\!\!\left(\rule[-1.5mm]{0mm}{4.75mm}\mysmallR^{\mbox{\tiny Fib$(n,k)$}}_{\mbox{\tiny ribbon}}\right)$ is edge-minimal and solitary. 
For the remainder of the proof, we focus on demonstrating that, when $k=3$ and $n \geq 3$, the DC relations on $L^{\mbox{\tiny Fib}}(n,k)$ uniquely determine edge products; Lemma 6.1.1 of \cite{DD} then forces all these edge products to be positive rational numbers. 
Thus, we can apply \SolitaryLatticeProposition\ to conclude that each $L^{\mbox{\tiny Fib}}(n,3)$ is edge-minimal. 

\begin{figure}[t]
\begin{center}
{\bf \FibSpineFig:}  {The spine and vertebra of $L^{\mbox{\tiny Fib}}(9,3)$.\\   (The particular choice of $n=9$ here is simply for illustrative purposes.)}

\vspace*{-0.15in}
\setlength{\unitlength}{1.5cm}
\begin{picture}(4,6.5)
\put(2.1,6){\scriptsize $(1,10,19)$}
\put(3.35,5){\scriptsize $(1,10,20)$}
\put(3.35,4){\scriptsize $(1,11,20)$}
\put(3.35,3){\scriptsize $(1,15,20)$}
\put(3.35,2){\scriptsize $(1,16,20)$}
\put(3.35,1){\scriptsize $(1,17,20)$}
\put(0.05,5){\scriptsize $(2,10,19)$}
\put(0.05,4){\scriptsize $(2,11,19)$}
\put(0.05,3){\scriptsize $(2,15,19)$}
\put(0.05,2){\scriptsize $(2,16,19)$}
\put(0.05,1){\scriptsize $(2,17,19)$}
\put(2.1,5){\scriptsize $(1,11,19)$}
\put(2.35,4){\tiny $(2,10,20)$}
\put(2.35,3){\tiny $(2,11,20)$}
\put(2.35,2){\tiny $(2,14,20)$}
\put(2.35,1){\tiny $(2,16,20)$}
\put(1.2,4){\tiny $(1,15,19)$}
\put(1.2,3){\tiny $(1,16,19)$}
\put(1.2,2){\tiny $(1,17,19)$}
\put(2.35,-0.1){\scriptsize $(2,17,20)$}
\put(0,0){\qbezier(2,6)(2.625,5.5)(3.25,5)}
\put(0,0){\qbezier(2,5)(2.625,4.5)(3.25,4)}
\put(0,0){\qbezier(2,4)(2.625,3.5)(3.25,3)}
\put(0,0){\qbezier(2,3)(2.625,2.5)(3.25,2)}
\put(0,0){\qbezier(2,2)(2.625,1.5)(3.25,1)}
\put(0,0){\qbezier(1,5)(1.625,4.5)(2.25,4)}
\put(0,0){\qbezier(1,4)(1.625,3.5)(2.25,3)}
\put(0,0){\qbezier(1,3)(1.625,2.5)(2.25,2)}
\put(0,0){\qbezier(1,2)(1.625,1.5)(2.25,1)}
\put(0,0){\qbezier(1,1)(1.625,0.5)(2.25,0)}
\put(0,0){\qbezier(2,6)(1.5,5.5)(1,5)}
\put(0,0){\qbezier(2,5)(1.5,4.5)(1,4)}
\put(0,0){\qbezier(2,4)(1.5,3.5)(1,3)}
\put(0,0){\qbezier(2,3)(1.5,2.5)(1,2)}
\put(0,0){\qbezier(2,2)(1.5,1.5)(1,1)}
\put(0,0){\qbezier(3.25,5)(2.75,4.5)(2.25,4)}
\put(0,0){\qbezier(3.25,4)(2.75,3.5)(2.25,3)}
\put(0,0){\qbezier(3.25,3)(2.75,2.5)(2.25,2)}
\put(0,0){\qbezier(3.25,2)(2.75,1.5)(2.25,1)}
\put(0,0){\qbezier(3.25,1)(2.75,0.5)(2.25,0)}
\put(0,0){\qbezier(2.25,4)(2.25,3.5)(2.25,3)}
\put(0,0){\qbezier(2.25,2)(2.25,1)(2.25,0)}
\put(0,0){\qbezier(3.25,3)(3.25,2)(3.25,1)}
\put(0,0){\qbezier(3.25,5)(3.25,4.5)(3.25,4)}
\put(0,0){\qbezier(2,6)(2,5.5)(2,5)}
\put(0,0){\qbezier(2,4)(2,3)(2,2)}
\put(0,0){\qbezier(1,5)(1,4.5)(1,4)}
\put(0,0){\qbezier(1,3)(1,2)(1,1)}
\multiput(1,4)(0,-0.2){6}{\line(0,-1){0.1}}
\multiput(2,5)(0,-0.2){6}{\line(0,-1){0.1}}
\multiput(2.25,3)(0,-0.2){6}{\line(0,-1){0.1}}
\multiput(3.25,4)(0,-0.2){6}{\line(0,-1){0.1}}
%
\put(1,5){\VertexBlankFib}
\put(1,4){\VertexBlankFib}
\put(1,3){\VertexBlankFib}
\put(1,2){\VertexBlankFib}
\put(1,1){\VertexBlankFib}
\put(2,6){\VertexBlankFib}
\put(2,4){\VertexBlankFib}
\put(2,3){\VertexBlankFib}
\put(2,2){\VertexBlankFib}
\put(2,5){\VertexBlankFib}
\put(2.25,4){\VertexBlankFib}
\put(2.25,3){\VertexBlankFib}
\put(2.25,2){\VertexBlankFib}
\put(2.25,1){\VertexBlankFib}
\put(2.25,0){\VertexBlankFib}
\put(3.25,5){\VertexBlankFib}
\put(3.25,4){\VertexBlankFib}
\put(3.25,3){\VertexBlankFib}
\put(3.25,2){\VertexBlankFib}
\put(3.25,1){\VertexBlankFib}
\put(1,5.04){\NEEdgeLabelForLatticeI{{\em 1}}}
\put(3.15,5.01){\NWEdgeLabelForLatticeI{{\em 1}}}
\put(0.85,3.89){\NEEdgeLabelForLatticeI{{\em 1}}}
\put(3.35,3.87){\NWEdgeLabelForLatticeI{{\em 1}}}
\put(0.85,4.28){\NEEdgeLabelForLatticeI{{\em 1}}}
\put(3.35,4.22){\NWEdgeLabelForLatticeI{{\em 1}}}
\put(0.85,2.89){\NEEdgeLabelForLatticeI{{\em 1}}}
\put(3.35,2.87){\NWEdgeLabelForLatticeI{{\em 1}}}
\put(0.85,3.28){\NEEdgeLabelForLatticeI{{\em 1}}}
\put(3.35,3.22){\NWEdgeLabelForLatticeI{{\em 1}}}
\put(0.85,1.89){\NEEdgeLabelForLatticeI{{\em 1}}}
\put(3.35,1.87){\NWEdgeLabelForLatticeI{{\em 1}}}
\put(0.85,2.28){\NEEdgeLabelForLatticeI{{\em 1}}}
\put(3.35,2.22){\NWEdgeLabelForLatticeI{{\em 1}}}
\put(0.85,0.89){\NEEdgeLabelForLatticeI{{\em 1}}}
\put(3.35,0.87){\NWEdgeLabelForLatticeI{{\em 1}}}
\put(0.85,1.28){\NEEdgeLabelForLatticeI{{\em 1}}}
\put(3.35,1.22){\NWEdgeLabelForLatticeI{{\em 1}}}
\put(1.15,0.04){\NEEdgeLabelForLatticeI{{\em 1}}}
\put(3.15,0.01){\NWEdgeLabelForLatticeI{{\em 1}}}
\put(2,5.04){\VerticalEdgeLabelForLatticeI{{\em 8}}}
\put(1,4){\VerticalEdgeLabelForLatticeI{{\em 8}}}
\put(3.25,4){\VerticalEdgeLabelForLatticeI{{\em 8}}}
\put(2.25,3){\VerticalEdgeLabelForLatticeI{{\em 8}}}
\put(2,3){\VerticalEdgeLabelForLatticeI{{\em 3}}}
\put(2,2){\VerticalEdgeLabelForLatticeI{{\em 2}}}
\put(1,2){\VerticalEdgeLabelForLatticeI{{\em 3}}}
\put(1,1){\VerticalEdgeLabelForLatticeI{{\em 2}}}
\put(2.25,1){\VerticalEdgeLabelForLatticeI{{\em 3}}}
\put(2.25,0){\VerticalEdgeLabelForLatticeI{{\em 2}}}
\put(3.25,2){\VerticalEdgeLabelForLatticeI{{\em 3}}}
\put(3.25,1){\VerticalEdgeLabelForLatticeI{{\em 2}}}
\end{picture}
\end{center}

\

\vspace*{-0.4in}
\end{figure}

Now take $k=3$ and $n \geq 3$ and suppose we are given a set of edge products $\{\myqP_{\selt,\telt}\}$ for $L^{\mbox{\tiny Fib}}(n,3)$ such that the set of scalar pairs $\{(\myqX_{\telt,\selt},\myqY_{\selt,\telt}) := (\sqrt{\myqP_{\selt,\telt}},\sqrt{\myqP_{\selt,\telt}})\}$ satisfies the DC relations. 
Let $x$ denote the edge product on edge $(2,n+1,2n+1) \myarrow{1} (1,n+1,2n+1)$. 
The {\em spine} of $L^{\mbox{\tiny Fib}}(n,3)$ is the chain $(1,2n-1,2n+1) \myarrow{2} (1,2n-2,2n+1) \myarrow{3} \cdots \mylongarrow{n-2} (1,n+2,2n+1) \mylongarrow{n-1} (1,n+1,2n+1)$. (For an example, see \FibSpineFig.) 
For a spinal edge $(1,2n+1-i,2n+1) \myarrow{i} (1,2n-i,2n+1)$ with $3 \leq i \leq n-1$, the corresponding `cube-shaped' {\em vertebra} is the eight-element Boolean lattice determined by the interval $[(2,2n+1-i,2n+2),(1,2n-i,2n+1)]$. 
Within this vertebra, it is easy to see that the DC relations require that all color $i$ edges have edge product equal to unity, all color $1$ edges incident with one of the vertices in the edge $(2,2n+1-i,2n+1) \myarrow{i} (2,2n-i,2n+1)$ have edge product $x$, and all color $1$ edges incident with one of the vertices in the edge $(1,2n+1-i,2n+2) \myarrow{i} (1,2n-i,2n+2)$ have edge product $2-x$. 

Now consider the vertebra with spinal edge $(1,2n-1,2n+1) \myarrow{2} (1,2n-2,2n+1)$.  
To make it easier to reference edges in this vertebra, we let $\relt_{1} := (2,2n-1,2n+2)$, $\selt_{1} := (2,2n-1,2n+1)$, $\telt_{1} := (1,2n-1,2n+2)$, and $\uelt_{1} = (1,2n-1,2n+1)$, and we let $\relt_{2} := (2,2n-2,2n+2)$, $\selt_{2} := (2,2n-2,2n+1)$, $\telt_{2} := (1,2n-2,2n+2)$, and $\uelt_{2} = (1,2n-2,2n+1)$. 
To see that the DC relations require the edge products in the following table, one must consider the requirements imposed by DC relations along the `back wall' of the lattice whose topmost edges form the chain $(n,n+2,2n+1) \mylongarrow{n-1} (n-1,n+2,2n+1) \mylongarrow{n-1} (n-1,n+1,2n+1) \mylongarrow{n-2} \cdots \myarrow{2} (2,n+1,2n+1) \myarrow{1} (1,n+1,2n+1)$. 
\begin{center}
{\small \begin{tabular}{|c|c|c|c|}
\hline
$\myqP_{\relt_{2},\selt_{2}} = x$ & $\myqP_{\selt_{2},\uelt_{2}} = x$ & $\myqP_{\relt_{2},\telt_{2}} = 2-x$ & $\myqP_{\telt_{2},\uelt_{2}} = 2-x$\\
\hline
\rule[-2.2mm]{-0.05mm}{6.5mm} & $\myqP_{\selt_{1},\selt_{2}} = \frac{n-1}{n-2}$ & $\myqP_{\telt_{1},\telt_{2}} = \frac{(n-1)(2-x)}{3(n-1)-(n-2)x}$ & $\myqP_{\uelt_{1},\uelt_{2}} = 1$\\
\hline
\rule[-2.2mm]{-0.05mm}{6.5mm}$\myqP_{\relt_{1},\selt_{1}} = \frac{(n-1)+(n-2)x}{n-1}$ & $\myqP_{\selt_{1},\uelt_{1}} = \frac{(n-2)x}{n-1}$ & $\myqP_{\relt_{1},\telt_{1}} = \frac{(n-2)(x)(3(n-1)-(n-2)x)}{(n-1)((n-1)+(n-2)x)}$ & $\myqP_{\telt_{1},\uelt_{1}} = \frac{3(n-1)-(n-2)x}{n-1}$\\
\hline
\end{tabular}
}
\end{center}
We omitted $\myqP_{\relt_{1},\relt_{2}}$ in the preceding table because the DC relations force this product to equal both of the quantities $\frac{(n-1)^{2}x}{(n-2)((n-1)+(n-2)x)}$ and $\frac{(n-1)^{2}(2-x)^{2}((n-1)+(n-2)x)}{(n-2)(x)(3(n-1)-(n-2)x)^{2}}$. 
Set these two quantities equal to one another and solve for $x$ to see that $x$ must be one of $-\frac{n-1}{n-2}$, $0$, or $\frac{n-1}{n}$. 
The first two of these three possibilities for $x$ will give us division by zero in one of the above expressions for $\myqP_{\relt_{1},\relt_{2}}$. 
But, since we were given a set of edge products that satisfy the DC relations, then these two possibilities for $x$ are ruled out. 
So, we must have $x = \frac{n-1}{n}$.

To complete the argument, we point out that just as we uniquely determined edge products for vertebra by working from the top of the spine to the bottom, edge products can be uniquely determined on all other cube-shaped and diamond-shaped intervals in this lattice.\hfill\QED

Enumerating ballot-admissible skew-shaped tableaux is known to be a difficult problem in general. 
However, for special families of shapes, the problem can be more tractable. 
This appears to be the case for the set $\mathcal{B}\!\left(\rule[-1.5mm]{0mm}{4.75mm}\mysmallR^{\mbox{\tiny Fib$(3,k)$}}_{\mbox{\tiny ribbon}}\right)$ of ballot-admissible tableaux with a Fibonacci ribbon shape. 
For $k \in \{1,2,\ldots,8\}$, we have empirically determined that the sets $\mathcal{B}\!\left(\rule[-1.5mm]{0mm}{4.75mm}\mysmallR^{\mbox{\tiny Fib$(3,k)$}}_{\mbox{\tiny ribbon}}\right)$ have the following sizes: 
\begin{center}
\begin{tabular}{|c||c|c|c|c|c|c|c|c|}
\hline
$k$ & $1$ & $2$ & $3$ & $4$ & $5$ & $6$ & $7$ & $8$\\
\hline
\rule[-2.25mm]{0mm}{6.5mm}$\rule[-1.5mm]{0.2mm}{4.95mm}\, \mathcal{B}\!\left(\rule[-1.5mm]{0mm}{4.75mm}\mysmallR^{\mbox{\tiny Fib$(3,k)$}}_{\mbox{\tiny ribbon}}\right)\rule[-1.5mm]{0.2mm}{4.95mm}$ & $1$ & $1$ & $2$ & $4$ & $8$ & $17$ & $37$ & $82$\\
\hline
\end{tabular}
\end{center}

We now consider certain sets of walks in $\mathbb{Z}^{2}$ from the origin $(0,0)$ to a point $(k,0)$ on the $x$-axis where each step of the walk adds some integer pair to the current position and the allowed steps are from some fixed set of vectors.\footnote{Such walks are called `bridges' or `excursions' in \cite{Betal}.} 
Call such a walk a {\em Motzkin path} if each step is one of $U$, $D$, or $L$, where $U$ is an upward diagonal that adds $(1,1)$ to the current position, $D$ is a downward diagonal adding $(1,-1)$, and $L$ is a level step adding $(1,0)$. 
Each such Motzkin path can be identified uniquely as a word of $k$ letters from the letter set $\{U,D,L\}$. 
Such a Motzkin path is {\em peakless} if an upward diagonal is never immediately followed by a downward diagonal, i.e.\ the corresponding word has no $UD$'s, and is {\em topside} if the path never descends below the $x$-axis. 
Let $\mathcal{M}(k)$ be the set of topside peakless Motzkin paths from $(0,0)$ to $(k,0)$. 
\FibMotzkinFig\ presents $\mathcal{M}(5)$, the set of eight topside peakless Motzkin paths from $(0,0)$ to $(5,0)$. 
It is well known that these sets of paths are enumerated by certain so-called generalized Catalan numbers, see OEIS-A004148. 
Since the `Motzkin numbers' enumerate all topside Motzkin paths (see OEIS-A001006), perhaps the numbers $\rule[-1.5mm]{0.2mm}{4.95mm}\, \mathcal{M}(k)\rule[-1.5mm]{0.2mm}{4.95mm}$ might be called `peakless Motzkin numbers', although Barry in \cite{PB} and Nkwanta in \cite{AN} suggest the  nomenclature of `RNA sequence' or `RNA numbers', due to Waterman's discovery of this sequence in enumerating secondary structures related to RNA \cite{Waterman}.
Here is a table of these numbers for small values of $k$: 
\begin{center}
\begin{tabular}{|c||c|c|c|c|c|c|c|c|}
\hline
$k$ & $1$ & $2$ & $3$ & $4$ & $5$ & $6$ & $7$ & $8$\\
\hline
\rule[-2.25mm]{0mm}{6.5mm}$\rule[-1.5mm]{0.2mm}{4.95mm}\, \mathcal{M}(k)\rule[-1.5mm]{0.2mm}{4.95mm}$ & $1$ & $1$ & $2$ & $4$ & $8$ & $17$ & $37$ & $82$\\
\hline
\end{tabular}
\end{center}

\begin{figure}[ht]
\begin{center}
{\bf \FibMotzkinFig:}  {\small Below are the eight top-side peakless Motzkin paths from $(0,0)$ to $(5,0)$.}

\setlength{\unitlength}{0.7cm}
\begin{picture}(24,2)
\thicklines
\multiput(0,0)(3,0){8}{
\thicklines
\put(0,0.5){\line(1,0){2.5}}
\put(0,0.5){\line(0,1){1}}
\thinlines
\multiput(0.5,0.5)(0.5,0){5}{\color{gray}\line(0,1){1}}
\multiput(0,1)(0,0.5){2}{\color{gray}\line(1,0){2.5}}
}
\thicklines
\put(0,0.5){\color{Red}\line(1,0){0.5}}
\put(0.5,0.5){\color{Red}\line(1,0){0.5}}
\put(1,0.5){\color{Red}\line(1,0){0.5}}
\put(1.5,0.5){\color{Red}\line(1,0){0.5}}
\put(2,0.5){\color{Red}\line(1,0){0.5}}
\put(0.1,0){\small \color{Red}$L$} 
\put(0.6,0){\small \color{Red}$L$} 
\put(1.1,0){\small \color{Red}$L$} 
\put(1.6,0){\small \color{Red}$L$} 
\put(2.1,0){\small \color{Red}$L$} 
\put(3,0.5){\color{Blue}\qbezier(0,0)(0.25,0.25)(0.5,0.5)}
\put(3.5,1){\color{Blue}\line(1,0){0.5}}
\put(4,1){\color{Blue}\qbezier(0,0)(0.25,-0.25)(0.5,-0.5)}
\put(4.5,0.5){\color{Blue}\line(1,0){0.5}}
\put(5,0.5){\color{Blue}\line(1,0){0.5}}
\put(3.1,0){\small \color{Blue}$U$} 
\put(3.6,0){\small \color{Blue}$L$} 
\put(4.1,0){\small \color{Blue}$D$} 
\put(4.6,0){\small \color{Blue}$L$} 
\put(5.1,0){\small \color{Blue}$L$} 
\put(6,0.5){\color{Green}\line(1,0){0.5}}
\put(6.5,0.5){\color{Green}\qbezier(0,0)(0.25,0.25)(0.5,0.5)}
\put(7,1){\color{Green}\line(1,0){0.5}}
\put(7.5,1){\color{Green}\qbezier(0,0)(0.25,-0.25)(0.5,-0.5)}
\put(8,0.5){\color{Green}\line(1,0){0.5}}
\put(6.1,0){\small \color{Green}$L$} 
\put(6.6,0){\small \color{Green}$U$} 
\put(7.1,0){\small \color{Green}$L$} 
\put(7.6,0){\small \color{Green}$D$} 
\put(8.1,0){\small \color{Green}$L$} 
\put(9,0.5){\color{Purple}\line(1,0){0.5}}
\put(9.5,0.5){\color{Purple}\line(1,0){0.5}}
\put(10,0.5){\color{Purple}\qbezier(0,0)(0.25,0.25)(0.5,0.5)}
\put(10.5,1){\color{Purple}\line(1,0){0.5}}
\put(11,1){\color{Purple}\qbezier(0,0)(0.25,-0.25)(0.5,-0.5)}
\put(9.1,0){\small \color{Purple}$L$} 
\put(9.6,0){\small \color{Purple}$L$} 
\put(10.1,0){\small \color{Purple}$U$} 
\put(10.6,0){\small \color{Purple}$L$} 
\put(11.1,0){\small \color{Purple}$D$} 
\put(12,0.5){\color{CarnationPink}\qbezier(0,0)(0.25,0.25)(0.5,0.5)}
\put(12.5,1){\color{CarnationPink}\line(1,0){0.5}}
\put(13,1){\color{CarnationPink}\line(1,0){0.5}}
\put(13.5,1){\color{CarnationPink}\qbezier(0,0)(0.25,-0.25)(0.5,-0.5)}
\put(14,0.5){\color{CarnationPink}\line(1,0){0.5}}
\put(12.1,0){\small \color{CarnationPink}$U$} 
\put(12.6,0){\small \color{CarnationPink}$L$} 
\put(13.1,0){\small \color{CarnationPink}$L$} 
\put(13.6,0){\small \color{CarnationPink}$D$} 
\put(14.1,0){\small \color{CarnationPink}$L$} 
\put(15,0.5){\color{Orange}\line(1,0){0.5}}
\put(15.5,0.5){\color{Orange}\qbezier(0,0)(0.25,0.25)(0.5,0.5)}
\put(16,1){\color{Orange}\line(1,0){0.5}}
\put(16.5,1){\color{Orange}\line(1,0){0.5}}
\put(17,1){\color{Orange}\qbezier(0,0)(0.25,-0.25)(0.5,-0.5)}
\put(15.1,0){\small \color{Orange}$L$} 
\put(15.6,0){\small \color{Orange}$U$} 
\put(16.1,0){\small \color{Orange}$L$} 
\put(16.6,0){\small \color{Orange}$L$} 
\put(17.1,0){\small \color{Orange}$D$} 
\put(18,0.5){\color{Goldenrod}\qbezier(0,0)(0.25,0.25)(0.5,0.5)}
\put(18.5,1){\color{Goldenrod}\line(1,0){0.5}}
\put(19,1){\color{Goldenrod}\line(1,0){0.5}}
\put(19.5,1){\color{Goldenrod}\line(1,0){0.5}}
\put(20,1){\color{Goldenrod}\qbezier(0,0)(0.25,-0.25)(0.5,-0.5)}
\put(18.1,0){\small \color{Goldenrod}$U$} 
\put(18.6,0){\small \color{Goldenrod}$L$} 
\put(19.1,0){\small \color{Goldenrod}$L$} 
\put(19.6,0){\small \color{Goldenrod}$L$} 
\put(20.1,0){\small \color{Goldenrod}$D$} 
\put(21,0.5){\color{ProcessBlue}\qbezier(0,0)(0.25,0.25)(0.5,0.5)}
\put(21.5,1){\color{ProcessBlue}\qbezier(0,0)(0.25,0.25)(0.5,0.5)}
\put(22,1.5){\color{ProcessBlue}\line(1,0){0.5}}
\put(22.5,1.5){\color{ProcessBlue}\qbezier(0,0)(0.25,-0.25)(0.5,-0.5)}
\put(23,1){\color{ProcessBlue}\qbezier(0,0)(0.25,-0.25)(0.5,-0.5)}
\put(21.1,0){\small \color{ProcessBlue}$U$} 
\put(21.6,0){\small \color{ProcessBlue}$U$} 
\put(22.1,0){\small \color{ProcessBlue}$L$} 
\put(22.6,0){\small \color{ProcessBlue}$D$} 
\put(23.1,0){\small \color{ProcessBlue}$D$} 
\end{picture}
\end{center}
\end{figure}
We have confirmed by computer experiments the coincidence of the sizes of the sets $\mathcal{B}\!\left(\rule[-1.5mm]{0mm}{4.75mm}\mysmallR^{\mbox{\tiny Fib$(3,k)$}}_{\mbox{\tiny ribbon}}\right)$ and $\mathcal{M}(k)$ for values of $k$ up to $k=20$. 
Based on this evidence, we proffer this conjecture:

\noindent 
{\bf \MotzkinConjecture}\ \ {\sl For all positive integers $k$, we have} $\rule[-1.5mm]{0.2mm}{4.95mm}\, \mathcal{B}\!\left(\rule[-1.5mm]{0mm}{4.75mm}\mysmallR^{\mbox{\tiny Fib$(3,k)$}}_{\mbox{\tiny ribbon}}\right)\rule[-1.5mm]{0.2mm}{4.95mm} = \rule[-1.5mm]{0.2mm}{4.95mm}\, \mathcal{M}(k)\rule[-1.5mm]{0.2mm}{4.95mm}$.  {\sl That is, the sets of ballot-admissible tableaux with Fibonacci ribbon shape} $\mysmallR^{\mbox{\tiny Fib$(3,k)$}}_{\mbox{\tiny ribbon}}$ {\sl and topside peakless Motzkin paths from $(0,0)$ to $(k,0)$ are equinumerous.}

Moreover, we ask:

\noindent
{\bf \MotzkinOpen}\ \ (1) Find a bijective proof of the equality proposed in \MotzkinConjecture. 
(2) What can be said about the enumeration of $\rule[-1.5mm]{0.2mm}{4.95mm}\, \mathcal{B}\!\left(\rule[-1.5mm]{0mm}{4.75mm}\mysmallR^{\mbox{\tiny Fib$(n,k)$}}_{\mbox{\tiny ribbon}}\right)\rule[-1.5mm]{0.2mm}{4.95mm}$ for $n > 3$?

%
\renewcommand{\refname}{\Large \bf References}
\renewcommand{\baselinestretch}{1.1}

\end{document}